\documentclass[11pt]{article}

\usepackage{epsfig}
\usepackage{graphicx}
\usepackage{amssymb}
\usepackage{amsmath}
\usepackage{color}
\usepackage{psfrag}

\usepackage[francais,english]{babel}
\usepackage[latin1]{inputenc}

\newtheorem{remark}{Remark}

\newtheorem{prop}{Proposition}

\newtheorem{lem}{Lemma}

\def\proof{ {{\bf Proof.}}}
\def\endproof{\hfill$\diamondsuit$\par\medskip}

\newcommand{\dps}[1]  {\displaystyle{#1} }

\newcommand{\uu}[1]  {{\boldsymbol #1} }
\newcommand{\diverg}[1]  {{\rm div} #1}

\def\R{\mathbb{R}}

\def\tr{\mathrm{tr}} 

\def\x{\uu{x}}
\def\y{\uu{y}}
\def\u{\uu{u}}
\def\v{\uu{v}}
\def\w{\uu{w}}
\def\f{\uu{f}}
\def\g{\uu{g}}
\def\n{\uu{n}}
\def\m{\uu{m}}
\def\t{\uu{t}}
\def\ttau{\uu{\tau}}
\def\ssigma{\uu{\sigma}}
\newcommand{\Om}{\Omega }
\newcommand{\ale}{\hat {\mathcal A}_t }
\newcommand{\invale}{\hat {\mathcal A}_t^{-1} }

\newcommand{\HH}     { \mathbb H      }
\newcommand{\Id}     { \rm Id      }
\newcommand{\dt}     { \delta t }
\newcommand{\alen}{ {\mathcal A}_{n,n+1} }
\newcommand{\invalen}{ {\mathcal A}_{n,n+1}^{-1} }
\newcommand{\invalenm}{ {\mathcal A}_{n-1,n}^{-1} }

\newcommand{\Frac}[2]{\displaystyle{\frac{#1}{#2}}} 
\newcommand{\Int}[2]{\displaystyle{\int_{#1}^{#2}}}

\newcommand{\comment}[1]{ }

\title{Generalized Navier Boundary Condition and Geometric Conservation Law for surface tension}
\author{J.-F. Gerbeau$^{1}$ and T. Lelièvre$^{2}$\\
 $^{1}$ INRIA,
 Rocquencourt, B.P.105,\\
78153 Le Chesnay Cedex, France \\
jean-frederic.gerbeau@inria.fr\\
 $^{2}$ CERMICS,
Ecole Nationale des Ponts et Chauss{\'e}es,\\
6 \& 8 Av. Pascal, 77455 Champs-sur-Marne, France\\
lelievre@cermics.enpc.fr}

\begin{document}

\selectlanguage{english}

\maketitle

\section*{Abstract}
We consider two-fluid flow problems in an Arbitrary Lagrangian
Eulerian (ALE) framework. The purpose of this work is twofold. First,
we address the problem of the moving contact line, namely the line
common to the two fluids and the wall. Second, we perform a stability
analysis in the energy norm for various numerical schemes, taking into
account the gravity and surface tension effects.

The problem of the moving contact line is treated with the so-called
Generalized Navier Boundary Conditions. Owing to these
boundary conditions, it is possible to circumvent the incompatibility
between the classical no-slip boundary condition and the fact that the
contact line of the interface on the wall is actually moving. 

The energy stability analysis is based in particular on an extension
of the Geometry Conservation Law (GCL) concept to the case of moving
surfaces. This extension is useful to study the contribution of the
surface tension. 

The theoretical and computational results presented in this paper allow
us to propose a strategy which offers a good compromise between
efficiency, stability and artificial diffusion.

\section{Introduction}

A difficult problem in the modelling of two-fluid flows in a bounded
domain concerns the displacement of the contact line, namely the
points which are at the intersection of the solid boundary of the domain and
the interface separating the two fluids. The difficulty comes from the
fact that:
\begin{itemize}
\item the interface follows the fluid motion: the normal velocity of a
  point on the interface is the normal velocity of the fluid particle
  at the same point,
\item the fluid particles near the boundary of the domain tend to have the same velocity
  as the points of the boundary.
\end{itemize}
Thus, if the velocity of the points on the boundary of the domain is
zero (classical no-slip condition for a viscous fluid on a fixed
wall), the moving contact line does not move: this is the so-called
moving contact line problem.  We refer to the review
papers~\cite{shikhmurzaev-97,qian-wang-sheng-06} for an introduction to this problem.

We have to keep in mind that the no-slip boundary condition for
viscous flows is only an \emph{approximation} of the Navier boundary
conditions (which relates the slip velocity relative to the wall to
the shear rate boundary). The Navier boundary conditions has been
assessed by molecular dynamics (MD) simulations, but fails in the
vicinity of the contact line. One might think that continuum
mechanics models (at the macroscopic level) are not able to represent the movement of the contact
line. Nevertheless, hybrid experiments reported in~\cite{hadjiconstantinou-99,hadjiconstantinou-99-JCP,thompson-robbins-93} which consists in
imposing the slip velocity obtained by MD simulation in a continuum
model, have shown that continuum models may provide good
results. Recently, Qian, Wang and Sheng
\cite{qian-wang-sheng-03,qian-wang-sheng-06} have proposed a boundary
conditions, the so-called Generalized Navier Boundary Condition
(GNBC), which is completly defined at the continuum level and which
seems to be in good agreement with MD results, including in the
vicinity of the contact line. We also refer to~\cite{ren-e-07} for a
careful study using MD, also showing that continuum models can be used
to described the dynamics of the contact line.

In~\cite{qian-wang-sheng-03,qian-wang-sheng-06}, this boundary
condition is used with a finite difference scheme and the displacement
of the interface between the two fluids is handled with a phase-field
approach. In this paper, we show that the GNBC is very natural with a
discretization method based on a variational formulation, like the
finite element method.  In addition, the ALE formulation we adopt is
convenient to study the energy balance and it uses less numerical
parameters than in a phase-field formulation (for which a free energy
for the order parameter needs to be introduced, for example). The main
drawback of the ALE method compared to other methods to follow a
moving interface (like volume of fluid
methods~\cite{hirt-nichols-81,gueyffier-li-nadim-scardovelli-zaleski-99},
level set methods~\cite{sussman-smereka-osher-94,sethian-smereka-03}
or phase-field formulation) is that it does not allow very large
motion of the interface without remeshing, and it does not allow a
change of topology of the domains occupied by each fluid. On the other
hand, it is generally admitted that it is the method of choice when a
precise computation of the position of the interface is required, with
good robustness and conservation properties.  The ALE formulation has
been used in number of applications involving two-fluid flows. We
refer for example to~\cite{gerbeau-le-bris-lelievre-book} for an
application to the modelling of aluminium electrolysis cells, and to
Section~\ref{sec:num} for an application to flows in narrow channels.

For efficiency and simplicity, many numerical schemes for two-fluid
flows decouple the fluid resolution and the movement of the
interface. The interface displacement is therefore solved explicitly with respect
to the computation of the fluid velocity and pressure. The variational formulation proposed in this study being
well-suited to study the discrete energy of the system, we propose to
quantify the spurious energy due the gravity and the surface tension
when the interface displacement is handled explicitly. The Geometric Conservation
Law (GCL) is a natural tool to perform this study. To analyse the
contribution of the surface tension terms, we will have to extend the
GCL to surface integral.

Here is the outline of the paper. In Section~\ref{sec:GNBC}, the GNBC
is introduced. The ALE formulation used to implement this boundary
condition is detailed in Section~\ref{sec:ALE}. We prove some lemmata
on Geometric Conservation Law in Section \ref{sec:surf-gcl}. We then
study the energy conservation properties of the numerical scheme in
Section~\ref{sec:energy}. Finally, Section~\ref{sec:num} is devoted to
some numerical experiments which illustrate the theoretical results
and investigate some alternative schemes.

\section{The Generalized Navier Boundary Condition}\label{sec:GNBC}

\begin{figure}[htbp]
\psfrag{Omega1}{$\Omega_1$}
\psfrag{Omega2}{$\Omega_2$}
\psfrag{Sigma}{$\Sigma$}
\psfrag{dOmega}{$\partial\Omega$}
\psfrag{no}{$\n_{\partial \Omega}$}
\psfrag{ns}{$\n_{\Sigma}$}
\psfrag{m}{$\m$}
\psfrag{ts}{$\t_{\partial \Sigma}$}
\psfrag{to}{$\t_{\partial \Omega}$}
\psfrag{theta}{$\theta$}
\centerline{\epsfig{file=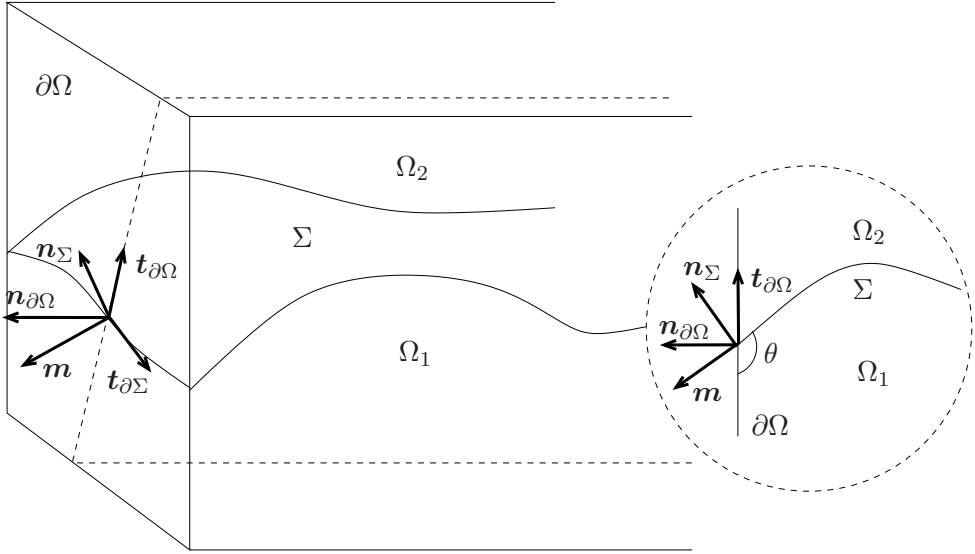,width=13cm}}
\caption{The domain $\Omega$ and various unit vectors. The 2d
  picture (in the dotted line circle) represents the vectors contained in the
  plane orthogonal to $\t_{\partial \Sigma}$.}\label{fig:normales}
\end{figure}

We are interested in the two-fluid Navier-Stokes equations, posed on a
bounded smooth domain $\Omega \subset \R^d$ (with $d=2$ or $d=3$), and
the time interval $(0,T)$:
\begin{equation}\label{eq:NS}
\left\{
\begin{array}{l}
\dps{\frac{\partial (\rho \u)}{\partial t}+\diverg(\rho \u \otimes \u) -
\diverg\left(\eta \left(\nabla \u + \nabla \u^T\right)\right)=-\nabla p
+ \gamma H \n_{\Sigma} \delta_{\Sigma} + \rho \g},\\
\diverg(\u)=0,\\
\dps{\frac{\partial \rho}{\partial t} + \diverg(\rho \u)=0}.
\end{array}
\right.
\end{equation}
The equation is posed in the distributional sense.  The velocity is
denoted by $\u$, the density by $\rho$, the viscosity by $\eta$ and
the pressure by $p$. The vector $\g$ denotes the gravity acceleration:
$$\g=-g \uu{e}_3$$
where  $(\uu{e}_1,\uu{e}_2,\uu{e}_3)$ is an orthonormal basis of the physical space. 
The term $\gamma H \n_{\Sigma} \delta_{\Sigma}$ is the surface
tension term that we will describe in detail below.  The system is
complemented by initial conditions $(\u(t=0),\rho(t=0))$. We suppose
that $\rho(t=0)$ takes two different values $\rho_1$ and
$\rho_2$. This property is then conserved as time evolves for the
function $\rho$ so that each fluid is distinguished from the other by
its density. In the following we denote by
\begin{equation}\label{eq:omega_i}
\Omega_i(t)=\left\{ \x \in \R^d, \, \rho(t,\x)=\rho_i \right\}
\end{equation}
the domain occupied at time $t$ by the fluid $i$. We suppose that
$\Omega_i(t)$ is a smooth domain, and we denote by
\begin{equation}\label{eq:int}
\Sigma(t)=\partial \Omega_1(t) \cap \partial \Omega_2(t)
\end{equation}
the interface between the two liquids. The viscosity $\eta$ may depend on the
fluid, so that $\eta=\eta(\rho)$ in~(\ref{eq:NS}). We denote in the
following by
\begin{equation}\label{eq:sigma}
\ssigma=\eta \left(\nabla \u + \nabla \u^T\right)
\end{equation}
the viscous stress tensor. In the surface tension term $\gamma H
\n_{\Sigma} \delta_{\Sigma}$, $\gamma$ is the surface tension
coefficient between the two fluids (which is supposed to be constant in
the following), $\n_{\Sigma}$ is the unit outward vector normal to
$\Omega_1$ (see Figure~\ref{fig:normales}) and $H$ is the mean curvature
of the interface $\Sigma$ positively counted with respect to the normal
$\n_{\Sigma}$. The distribution $\delta_{\Sigma}$ is defined by: for any
smooth function $\psi$
\begin{equation}\label{eq:delta}
\langle  \delta_{\Sigma}, \psi \rangle = \int_{\Sigma} \psi d \sigma_{\Sigma}
\end{equation}
where $\sigma_{\Sigma}$ denotes the Lebesgue measure ({\em i.e.} the
surface measure) on $\Sigma$.

Let us now describe the boundary conditions. We first suppose the
non-penetration condition:
\begin{equation}\label{eq:BC_u_n}
(\u - \u^b) \cdot \n_{\partial \Omega}=0 \text{ on $\partial\Omega$},
\end{equation}
where $\n_{\partial \Omega}$ denotes the unit outward vector normal to
$\Omega$ (see Figure~\ref{fig:normales}) and $\u^b$ is the velocity of
the boundary ($\u-\u^b$ is the slip velocity, $\u^b=0$ for a fixed
wall). In addition, we will suppose for the sake of simplicity that 
$$
\u^b\cdot\n_{\partial \Omega}=0,
$$
which implies that the the fluid lives in a fixed domain $\Omega$.


The fluid being viscous, we have to set another boundary condition on
$\partial\Omega$. As explained in the introduction, the classical
no-slip boundary condition ($\u = \u^b$ on $\partial\Omega$) would stick the interface on the wall, which is
clearly unacceptable. 

It may be worth recalling that the no-slip boundary condition is in
fact an \emph{approximation} of the Navier boundary condition:
\begin{equation}
  \label{eq:navier}
  \beta \left( \u - \u^b\right)\cdot \ttau + \ssigma \n_{\partial \Omega} \cdot \ttau=0,
\end{equation}
where $\beta$ is the slip coefficient. The coefficient $\beta$ is in
practice very large, which explains that $(\u-\u^b).\ttau = 0$ is a
good approximation in many cases. 

As explained in the introduction, the Navier boundary condition is not
appropriate to describe the moving contact line. This is the
motivation of the Generalized Navier Boundary Conditions (GNBC)
introduced in \cite{qian-wang-sheng-03}. To state the GNBC, we need
to define (see Figure~\ref{fig:normales}) the following vectors,
defined on the boundary $\partial \Sigma$ of the interface:
$\t_{\partial \Sigma}=\n_{\Sigma} \times \n_{\partial \Omega}$ the
tangent vector to $\partial \Sigma$, $\m=\t_{\partial \Sigma} \times
\n_{\Sigma}$ and $\t_{\partial \Omega}=\n_{\partial \Omega} \times
\t_{\partial \Sigma}$. Both sets of vectors $(\t_{\partial
  \Sigma},\n_{\Sigma},\m)$ and $(\t_{\partial \Sigma},\t_{\partial
  \Omega},\n_{\partial \Omega})$ are positively oriented orthonormal
basis. The GNBC writes: for any vector $\ttau$ tangent to $\partial
\Omega$,
\begin{equation}\label{eq:GNBC}
\beta \left( \u - \u^b\right) \cdot \ttau + \ssigma \n_{\partial \Omega} \cdot \ttau + \gamma \left(
    \m \cdot \t_{\partial \Omega}  - \cos(\theta_s) \right)
  \t_{\partial \Omega} \cdot \ttau \, \delta_{\partial \Sigma}=0,
\end{equation}
where as above $\beta$ is the slip coefficient, $\gamma$ is the
surface tension coefficient between the two fluids, and $\theta_s$ is
the static contact angle at the solid surface. The distribution
$\delta_{\partial \Sigma}$ is defined accordingly
with~(\ref{eq:delta}) by: for any smooth function $\psi$
\begin{equation}\label{eq:delta_}
\langle \delta_{\partial \Sigma},  \psi\rangle = \int_{\partial \Sigma} \psi d
l_{\partial \Sigma}
\end{equation}
where $l_{\partial \Sigma}$ denotes the Lebesgue measure ({\em i.e.} the
length measure) on the curve $\partial \Sigma$.


Compared to the classical Navier boundary condition~\eqref{eq:navier}, we
notice the extra term
$$
\left( \m \cdot \t_{\partial \Omega} -
  \cos(\theta_s) \right),
$$
which is called the \textit{uncompensated Young stress}
in~\cite{qian-wang-sheng-06}. This quantity measures the difference
between the dynamic contact angle $\theta$ made by the interface
$\Sigma$ and the boundary $\partial \Omega$ (we choose the convention
that this angle is measured in the fluid $1$, see
Figure~\ref{fig:normales}) and the static contact angle $\theta_s$.
In general, the definition of $\theta_s$ is a difficult problem both
from the experimental and the theoretical standpoints. In our
framework, $\theta_s$ is assumed to be a part of the data. The
uncompensated Young stress is concentrated along the boundary of the
interface and vanishes when the dynamic contact angle is equal to the
static contact angle.


\section{Variational formulation and discretization}\label{sec:ALE}

The aim of this section is to derive a variational formulation of the
system of equations~(\ref{eq:NS}), together with the boundary conditions~(\ref{eq:BC_u_n})--(\ref{eq:GNBC}).

\subsection{The weak ALE formulation}

\begin{figure}[htbp]
\psfrag{A}{$\ale$}
\psfrag{hO2}{$\hat\Omega_2$}
\psfrag{hO1}{$\hat\Omega_1$}
\psfrag{hS}{$\hat\Sigma$}
\psfrag{O2}{$\Omega_{2,t}$}
\psfrag{O1}{$\Omega_{1,t}$}
\psfrag{S}{$\Sigma_t$}
  \centerline{
    \epsfig{figure=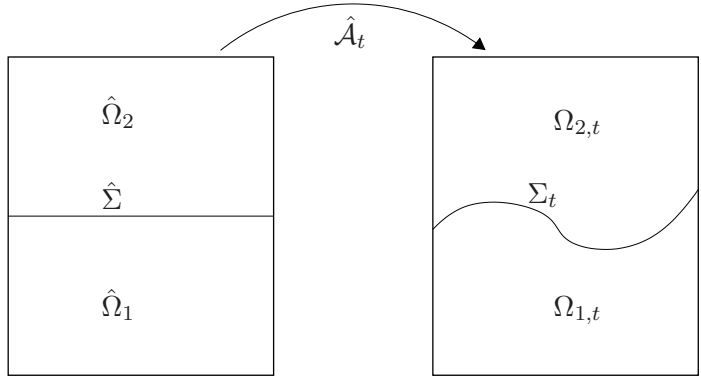,height=5cm,angle=0}}
    \caption{The partition of the domain $\Omega$ and the ALE mapping.}
    \label{fig:ale}
\end{figure}

We assume that for any time $t \geq 0$, there exists a smooth and
bijective mapping $\ale$ from a reference domain $\hat{\Omega}$
(divided into two separate subdomains $\hat{\Omega}_1$ and
$\hat{\Omega}_2$ such that
$\overline{\hat{\Omega}}=\overline{\hat{\Omega}_1}\cup\overline{\hat{\Omega}_2}$)
to the current domain $\Omega$ such that
$\ale(\hat{\Omega}_i)=\Omega_i(t)$ (see Figure~\ref{fig:ale}). The
inverse function (with respect to the space variable) of $\ale$ is
denoted $\invale$.

The velocity of the domain $\hat \w$ is defined by:
\begin{equation}\label{eq:w_hat}
\hat \w(t, \hat \x)=\frac{\partial}{\partial t} \ale(\hat \x).
\end{equation}
For any function $\psi(t,.)$ defined on $\Omega$, we denote by
$\hat{\psi}(t,.)$ the corresponding function defined on the reference domain
$\hat{\Omega}$ by
\begin{equation}\label{eq:psi_hat}
\hat{\psi}(t,\hat{\x})=\psi(t,\ale(\hat{\x})).
\end{equation}
For example, the velocity of the domain $\w$ on the current frame is
defined by
\begin{equation}\label{eq:w}
\w(t,\x)=\hat \w(t,\ale^{-1}(\x)).
\end{equation}
Notice that the functions $\psi$ and $\hat{\psi}$ are such that:
\begin{equation}\label{eq:chain_rule}
\frac{\partial \hat{\psi}}{\partial
  t}(t,\hat{\x})=\frac{\partial \psi}{\partial t} (t,\ale(\hat{\x}))+\w(t,\ale(\hat{\x}))\cdot \nabla \psi(t,\ale(\hat{\x})).
\end{equation}
The fact that $\ale$ maps $\hat{\Omega}_i$ to $\Omega_{i}(t)$ ($i=1$ or
$2$) implies that the velocity of the domain satisfies
\begin{equation}\label{eq:w_bord}
\w \cdot \n_i=\u \cdot \n_i \mbox{ on $\partial \Omega_i$, }
\end{equation}
where $i=1$ or $i=2$ and $\n_i$ denotes the unit outward vector normal to $\Omega_i$.
The density $\rho$ of the fluid is such that:
\begin{equation}\label{eq:weak-density_}
\rho(t,\x)=\hat \rho(\ale^{-1}(\x)),
\end{equation}
where $\hat \rho$ is equal to $\rho_1$ on $\hat \Omega_1$ and
$\rho_2$ on $\hat \Omega_2$.

The following functional spaces will be needed, respectively for the
velocity $\u$ and the pressure~$p$:
$$
V = L^2(0,T;\HH^1_\n(\Omega)), \quad M = L^2(0,T;L^2_0(\Omega)),
$$
where
$$
\HH^1_\n(\Omega) = \left\{\u \in \left(H^1(\Omega)\right)^d, \u \cdot \n_{\partial \Omega}=0 \mbox{ on $\partial \Omega$}\right\},
$$
and
$$L^2_0(\Omega)=\left\{ p \in L^2(\Omega), \, \int_\Omega p = 0\right\}.$$

We also introduce the test function spaces on the reference domain
$$
\hat V = \HH^1_\n (\hat \Omega), 
\qquad \hat M = L^2_0(\hat \Omega).
$$
In the moving frame, the test function spaces are defined by
$$
V_T = \{\v : [0,T] \times \Omega \rightarrow \R^d, \, \v(t,\x) = \hat \v(
\invale(\x)), \, \hat\v \in \hat V \},
$$
$$
M_T = \{q : [0,T] \times \Omega \rightarrow \R, \, q(t,\x) = \hat q(
\invale(\x)), \,\hat q \in \hat M \}.
$$
Thus, the test functions do not depend on time
in the reference frame $\hat\Omega$ whereas they do on the current
one: more precisely, let $\v$ be in $V_T$, then for a fixed $\hat \x
\in \hat\Omega$, $\v(t,\ale(\hat \x))$ does not depend on time while for a fixed $\x \in \Omega$, $\v(t,\x)$ does.

We are now in position to state the weak ALE formulation. It is the following coupled problem: we look for a
function $\ale:\hat\Omega \to \Omega$ and $(\u,p)$ in $V \times M$ such
that $\u(t=0,.)=\u_0$ and:
\begin{itemize}
\item The function $\ale$ is smooth and maps $\hat{\Omega}_i$ to
  $\Omega_{i}(t)$ ($i=1$ or $2$). The domains $\Omega_{i}(t)$ occupied by each
  fluid are thus defined by $\ale$ and the density of the fluid $\rho$ is defined by:
  \begin{equation}
    \label{eq:weak-density}
    \rho(t,\x) = \hat \rho(\ale^{-1}(\x)) = \rho_i, \qquad
    \mbox{ for }  \x \in \Omega_{i}(t).
  \end{equation}
\item For all $(\v,q)$ in $V_T \times M_T$,
  \begin{equation}\label{eq:weak-NS}
\left\{
\begin{array}{l}
    \dps{ \frac{d}{dt} \int_\Omega \rho \u \cdot \v + \int_\Omega \rho (\u-
    \w)\cdot \nabla \u \cdot \v  - \int_\Omega \diverg(\w) \rho \u \cdot
    \v }\\
 \quad \dps{+  \int_\Omega \frac{\eta}{2} \left( \nabla \u + \nabla \u^T \right):
    \left( \nabla \v + \nabla \v^T \right)
    - \int_\Omega p \, \diverg(\v) }  \\
\quad  \dps{=   -\gamma \int_{\Sigma} \tr(\nabla_\Sigma
  \v) \, d\sigma_{\Sigma} -   \beta \int_{\partial \Omega} (\u - \u^b)
  \cdot \v } \\
\qquad  \dps{+ \gamma
\int_{\partial \Sigma}  \cos(\theta_s)   \t_{\partial \Omega} \cdot \v
\, dl_{\partial \Sigma}+ \int_\Omega \f \cdot \v,}\\    
    \dps{\int_\Omega q \, \diverg(\u) = 0.}
  \end{array}
\right.
 \end{equation}
\end{itemize}
Of course, $\ale$ is not uniquely defined by the
condition~\eqref{eq:weak-density}: this is the arbitrary feature of the
ALE scheme. The function $\ale$ will be precisely defined at the
discrete level in Section~\ref{sec:complete_algo}.

{\bf Notation:} In~\eqref{eq:chain_rule},~\eqref{eq:weak-NS} and in the sequel, the
spatial differential operators are taken with respect to the Eulerian variable
$\x$. We omit to denote this explicitly for conciseness.

\subsection{Derivation of the weak ALE formulation}

Let us explain how this weak ALE formulation is obtained from the strong
formulation~(\ref{eq:NS}), with the boundary
conditions~(\ref{eq:BC_u_n})--(\ref{eq:GNBC}). 

The starting point of this derivation is based on the following Lemma.
\begin{lem}\label{lem:reynolds}
For any smooth function $\psi$ depending on time $t$ and space $\x$,
and any smooth function $\phi$ such that $\hat \phi$ (defined by $\hat{\phi}(t,\hat{\x})=\phi(t,\ale(\hat{\x}))$) is time-independent, we have:
\begin{eqnarray}
\lefteqn{\frac{d}{dt} \int_{\Omega} \psi(t,\x)\, \phi(t,\x)\,d\x}\label{eq:reynolds}\\
&=&\int_{\Omega} \phi(t,\x) \frac{\partial \psi}{\partial t}(t,\x)  + \phi(t,\x) \w(t,\x)\cdot \nabla \psi(t,\x) +\phi(t,\x) \diverg(\w(t,\x)) \psi(t,\x)\, d\x. \nonumber
\end{eqnarray}
\end{lem}

The first line in~(\ref{eq:weak-NS}) is obtained by multiplying the
material derivative in the equation on $\u$ in~(\ref{eq:NS}) by the test function $\v \in V_T$ and
integrating over $\Omega$:
\begin{align*}
\int_\Omega \frac{\partial (\rho \u)}{\partial t} \cdot \v &+\diverg(\rho
\u \otimes \u) \cdot \v =\int_\Omega \rho \frac{\partial \u}{\partial t} \cdot \v + \rho
\u \cdot \nabla \u \cdot \v,\\ 
&=\frac{d}{dt} \int_\Omega \rho \u \cdot \v - \int_\Omega \rho \w \cdot
\nabla \u \cdot \v - \diverg(\w) \rho \u \cdot \v + \rho \u \cdot \nabla \u \cdot \v,
\end{align*}
where we used successively the equation on $\rho$ in~(\ref{eq:NS}) and
then~(\ref{eq:reynolds}). The weak formulation of the terms involving
the pressure are classically obtained by integration by parts. It is straightforward to obtain the following variational
formulation for the term involving
the viscous stress and the surface tension term:
\begin{equation}\label{eq:1}
 \int_\Omega \frac{\eta}{2} \left(\nabla \u + \nabla
    \u^T\right) : \left( \nabla \v + \nabla
    \v^T\right)- \int_{\partial \Omega} \ssigma
  \n_{\partial \Omega} \cdot \v - \int_{\Sigma} \gamma H \v \cdot
  \n_\Sigma d \sigma_\Sigma.
\end{equation}

We now use the surface divergence formula (see~\cite{weatherburn-47} Equation~(24) p.~239 or~\cite{ambrosio-soner-96}, Equation~(3.8)). For any smooth hypersurface
$\Sigma$ in $\R^d$ (\emph{i.e.} a submanifold of $\R^d$ with codimension
$1$) with a smooth boundary $\partial \Sigma$ and normal $\n_\Sigma(\x)$ at point $\x$, one has: for any smooth
function $\uu{\Phi} : \Sigma \to \R^d$,
\begin{equation}\label{eq:div_th}
- \int_{\Sigma} H \uu{\Phi}\cdot \n_\Sigma \, d\sigma_{\Sigma}=\int_{\Sigma}\tr(\nabla_\Sigma \uu{\Phi}) \, d\sigma_{\Sigma} - \int_{\partial \Sigma} \uu{\Phi}\cdot \m
\, dl_{\partial \Sigma},
\end{equation} 
where the surface gradient $\nabla_\Sigma$ is defined by: for any smooth
vector field $\uu{X}$,
\begin{equation}\label{eq:surf_gradient}
\nabla_\Sigma \uu{X}=P_\Sigma(\x) \nabla \uu{X},
\end{equation} 
where $P_\Sigma(\x)$ is the orthogonal projector onto the tangent space
to $\Sigma$ at point $\x$:
$$P_\Sigma(\x)=\Id - \n_\Sigma(\x) \otimes \n_\Sigma(\x).$$
Notice that the surface gradient of $\uu{X}$ only depends on the values
of $\uu{X}$ on the surface $\Sigma$. The vector $\m$ is the normal vector to $\partial \Sigma$ in the
tangent space of $\Sigma$ pointing outwards of~$\Sigma$ (see Figure~\ref{fig:normales}).  The measure
$l_{\partial_\Sigma}$ is the Lebesgue measure on $\partial \Sigma$.

Using the surface divergence formula~(\ref{eq:div_th}), the last two
terms in~(\ref{eq:1}) writes:
\begin{align}
- \int_{\partial \Omega} \ssigma
  \n_{\partial \Omega} \cdot \v - \int_{\Sigma} \gamma H \v \cdot
  \n_\Sigma d \sigma_\Sigma
&=- \int_{\partial \Omega} \ssigma
  \n_{\partial \Omega} \cdot \v + \gamma \int_{\Sigma} \tr(\nabla_\Sigma
  \v) \, d\sigma_{\Sigma} \nonumber \\
&\quad - \gamma \int_{\partial \Sigma} \v \cdot \m
\, dl_{\partial \Sigma} \label{eq:2}.
\end{align}
We now use the GNBC~(\ref{eq:GNBC}) to rewrite the first and last terms
in the right-hand side of~(\ref{eq:2}):
\begin{align*}
- &\int_{\partial \Omega} \ssigma
  \n_{\partial \Omega} \cdot \v - \gamma \int_{\partial \Sigma} \v \cdot \m
\, dl_{\partial \Sigma}\\
&=  \beta \int_{\partial \Omega} (\u - \u^b) \cdot \v + \gamma
\int_{\partial \Sigma} \left(\m \cdot\t_{\partial \Omega} - \cos(\theta_s)  \right) \t_{\partial \Omega} \cdot \v
\, dl_{\partial \Sigma}  - \gamma \int_{\partial \Sigma} \v \cdot \m
\, dl_{\partial \Sigma},\\
&=  \beta \int_{\partial \Omega} (\u - \u^b) \cdot \v - \gamma
\int_{\partial \Sigma}  \cos(\theta_s)   \t_{\partial \Omega} \cdot \v
\, dl_{\partial \Sigma},
\end{align*}
where we have used the fact that $\v \cdot \m=\left(\v \cdot \t_{\partial
  \Omega}\right) \left( \m \cdot \t_{\partial \Omega} \right)$ (since $\v \cdot
   \n_{\partial \Omega}=0$).



With the divergence formula, we eliminate the mean curvature $H$ (which
is difficult to approximate at the discrete level), and we naturally
enforce the GNBC.

\subsection{Discretization}

The discretization is based on a finite element method in space, and an implicit Euler
time-discretization. The domain $\overline{\Omega}^n=\overline{\Omega}^n_1 \cup \overline{\Omega}^n_2$ at
the beginning of the $n$-th timestep, where $\Omega^n_i$ is the domain
occupied by the fluid $i$ at time $t_n$, plays the role of the reference
domain $\overline{\hat \Omega}=\overline{\hat\Omega}_1 \cup
\overline{\hat\Omega}_2$.

Given the mesh ${\mathcal M}^n={\mathcal M}_1^n \cup {\mathcal M}_2^n$ of
the domain\footnote{and therefore the domains occupied by each fluid}
$\overline{\Omega}^n=\overline{\Omega}_1^n \cup \overline{\Omega}_2^n$
and the
velocity $\u^n$ discretized in a finite element space at time $t_n$, we aim to
propagate these two items to time $t_{n+1}$, using the weak ALE
formulation~(\ref{eq:weak-NS}).

In addition to $({\mathcal M}^n,\u^n)$, let us give ourselves a space
discretization of the domain velocity $\w^n$ at time $t_n$.  Depending
on the scheme (implicit or explicit treatment of the displacement of the domain) the
velocity $\w^n$ may depend on $\u^n$ or on $\u^{n+1}$. We will come
back to its computation below, in Section~\ref{sec:complete_algo}. We
introduce the application
\begin{equation}\label{eq:alen}
\alen:\left\{
\begin{array}{ccc}
(\Omega^n_i)_{i=1,2} & \to & (\Omega^{n+1}_i)_{i=1,2} \\
\y & \mapsto &\x = \y + \delta t \,\w^{n}(\y)
\end{array}
\right.
,
\end{equation}
which might be seen as an approximation of $\hat {\mathcal A}_{t_{n+1}}
\circ \hat {\mathcal A}_{t_{n}}^{-1}$. This application defines the
domain occupied by each fluid\footnote{$\alen$ also defines
  the mesh at time $t_{n+1}$: for $i=1,2$, each node of ${\mathcal M}^n_i$ is
transported from $\Omega^n_i$ to $\Omega^{n+1}_i$ by $\alen$, thus defining
the mesh ${\mathcal M}^{n+1}_i$ of $\Omega^{n+1}_i$ at time $t_{n+1}$.} at time $t_{n+1}$:
$\Omega^{n+1}_i=\alen(\Omega^n_i)$, for $i=1,2$. Without loss of
generality, the time-step $\delta t=t_{n+1}-t_n$ is supposed to be constant.
In the sequel, our convention is that $\y$ denotes a point
in $(\Omega^n_i)_{i=1,2}$ and $\x$ a point in
$(\Omega^{n+1}_i)_{i=1,2}$.

\subsubsection{Discretization in space}

We consider a finite element discretization of the domain
$(\Omega^n_i)_{i=1,2}$. It is transported by the application $\alen$ to
a finite element discretization of the domain
$(\Omega^{n+1}_i)_{i=1,2}$. The finite element spaces at time $t_n$ for the velocity
and the pressure are respectively denoted by
$$
V_{h,n} \subset \HH^1_\n(\Omega), \qquad \qquad M_{h,n} \subset
L^2_0(\Omega).
$$
These finite element spaces depend on the time index $n$, since the
mesh is moving. They are supposed to satisfy the inf-sup condition
(see~\cite{ern-guermond-04,girault-raviart-86} for example).

As explained above, we use test functions which follow the
deformation of the domain given by $\alen$: the test functions at time $t_{n+1}$
belong to the following spaces:
$$
V_{h,n+1} = \{\v(t_{n+1},.): \Omega \rightarrow \R^N, \v(t_{n+1},\x) = \v(t_n,
\invalen(\x)), \v(t_n,.) \in V_{h,n} \},
$$
$$
M_{h,n+1} = \{q(t_{n+1},.): \Omega \rightarrow \R, q(t_{n+1},\x) =  q(t_n,
\invalen(\x)), q(t_n,.) \in  M_{h,n} \}.
$$
Unless there is a risk of confusion, we omit the index $h$ for the
functions belonging to the finite element spaces
$V_{h,n}$ or $M_{h,n}$.

\subsubsection{Time discretization and linearization}
We use the following semi-implicit Euler discretization of~(\ref{eq:weak-NS}): for a given
$\u^n \in V_{h,n}$, $(\Omega^n_i)_{i=1,2}$, $\w^n$ and
$(\Omega^{n+1}_i)_{i=1,2}$, compute $(\u^{n+1},p^{n+1}) \in V_{h,n+1}
\times M_{h,n+1}$ such that, for all $(\v(t_n,.),q(t_n,.)) \in V_{h,n}
\times M_{h,n}$,
\begin{equation}\label{eq:weak-NS-disc}
  \left\{
    \begin{array}{l}
      \dps{\frac{1}{\delta t} \int_{\Omega^{n+1}} \rho \u^{n+1} \cdot \v +
        \int_{\Omega^{n+1}} \rho (\u^n - \w^n) \cdot \nabla \u^{n+1} \cdot \v -
        \int_{\Omega^{n+1}} \diverg(\w^n) \rho \u^{n+1} \cdot \v} \\
      \dps{\quad+ \int_{\Omega^{n+1}} \frac{\eta}{2}(\nabla \u^{n+1} + (\nabla
        \u^{n+1})^T):(\nabla \v + \nabla \v^T) - \int_{\Omega^{n+1}} p^{n+1}
        \diverg(\v)}\\
      \dps{\quad + \int_{\Omega^{n+1}} \frac{\rho}{2} \,   \diverg(\u^n) \, \u^{n+1} \cdot
        \v + \frac{\delta \rho}{2}\int_{\Sigma^{n+1}} (\u^n-\w^n)\cdot
        \n_{\Sigma^{n+1}} \, \u^{n+1} \cdot \v \, d\sigma_{\Sigma^{n+1}} } \\
      \dps{\qquad = \frac{1}{\delta t} \int_{\Omega^{n}} \rho \u^n \cdot \v -
        \gamma \int_{\Sigma^{n+1}} \tr(\nabla_{\Sigma^{n+1}} \v) d\sigma_{\Sigma^{n+1}} - \beta
        \int_{\partial \Omega^{n+1}} (\u^{n+1}-\u^b) \cdot \v}\\
      \dps{\qquad \quad + \gamma \int_{\partial \Sigma^{n+1}} \cos(\theta_s)
        \t_{\partial \Omega} \cdot \v \, dl_{\partial \Sigma^{n+1}} +
        \int_{\Omega^{\star}} \rho \g \cdot
        \v,}\\
      \dps{\int_{\Omega^{n+1}} q \diverg(\u^{n+1})=0.}
    \end{array}
  \right.
\end{equation}
Notice that the superscript $n$ (in $\Omega^n$) emphasizes that we
consider the domain at time $t^n$, even if the boundary of the domain
is \emph{not} moving. When we integrate over $\Omega^{n+1}$, this is
to indicate that the test functions and the functions $\rho$, $\eta$
(whose values are deduced from the domains occupied by each fluid) are
taken at time~$t_{n+1}$. If a function defined on $\Omega^n$ appears
in an integral over $\Omega^{n+1}$, it means that this function is
transported on $\Omega^{n+1}$ by~$\alen$. For example,
$$
\int_{\Omega^{n+1}} \frac{\rho}{2} \,   \diverg(\u^n) \, \u^{n+1} \cdot
\v=\sum_{i=1}^2 \frac{\rho_i}{2} \int_{\Omega_i^{n+1}} \diverg(\u^n
\circ \invalen) \, \u^{n+1} \cdot \v(t_{n+1},.).
$$

The discretization~(\ref{eq:weak-NS-disc}) is obtained from the weak
ALE formulation~(\ref{eq:weak-NS}). In the third line
of~(\ref{eq:weak-NS-disc}) appear two terms which are required for
better stability properties of the scheme. The term
$\dps{\int_{\Omega^{n+1}} \frac{\rho}{2} \, \diverg(\u^n) \,
  \u^{n+1} \cdot \v}$ is standard. It is analogous to the well-known
modification introduced by Temam of the convective term (see
Section~III.5 in~\cite{temam-79}) which allows to recover at the
discrete level the skew-symmetry property of the advection term. The
second term $\dps{\frac{\delta \rho}{2}\int_{\Sigma^{n+1}}
  (\u^n-\w^n)\cdot \n_{\Sigma} \, \u^{n+1} \cdot \v \, d\sigma}$
where we have used the notation

\begin{equation}\label{eq:delta_rho}
  \delta \rho=\rho_2 - \rho_1,
\end{equation}
is in the same vein, but is specific to the context of two-fluid
flows. Notice that both these terms are strongly consistent: they
vanish for the exact solution.  They are introduced in order to
reproduce at the discrete level the energy estimates that can be
derived at the continuous level (see Section~\ref{sec:energy}).

The body force is integrated on a domain denoted by $\Om^\star$. We
will consider in the sequel different choices: $\Om^\star = \Om^n$ or
$\Om^{n+1}$ or $\Om^{n+1/2}$. By convention, we \emph{define} the
quantity $\int_{\Omega^{n+1/2}} \rho \g \cdot\v$ as $$\int_{\Omega^{n+1/2}}
\rho \g \cdot\v=
\frac{1}{2}
\int_{\Omega^{n+1}} \rho \g \cdot\v + \frac{1}{2} \int_{\Omega^{n}}
\rho \g \cdot\v.$$ It is worth noticing that the three choices have the same
computational cost as soon as the mesh velocity $\w^n$ is explicitly
obtained from the known velocity $\u^n$.They nevertheless lead to very
different behaviors, as shown in Section~\ref{sec:num-energy}.

\subsubsection{The complete algorithm}\label{sec:complete_algo}

To complete the presentation of the numerical scheme, it remains to
describe how the domain velocity~$\w^n$ is computed. The basic
requirement is the kinematic condition~(\ref{eq:w_bord}), which
ensures in particular that the nodes of the mesh which are initially
on the interface remain on the interface. In addition, $\alen$ defined
from $\w^n$ by~(\ref{eq:alen}) must be sufficiently smooth so that the
mesh remains regular enough for finite element computations.

In the practical problems we are interested in, it is sufficient to
adopt the very standard method that consists in solving a simple
Poisson problem to compute the velocity of the mesh.  Moreover, we
choose the displacement to be along one direction only, so that we
actually solve a \emph{scalar} Poisson problem.  This choice, which is
definitely reasonable in the physical situations that we consider, has
important favorable consequences on the quality of the algorithm.
This will be made precise in Section~\ref{sec:energy_disc}. In
addition, we discretize the velocity of the domain $\w^n$ in space
using the same finite element space as for the components of $\u^n$.

We may now write the complete algorithm. Suppose that
$(\Omega^n_i)_{i=1,2}$ and $(\u^n,p^n)$ are known. Then $\w^n$,
$(\Omega^{n+1}_i)_{i=1,2}$ and $(\u^{n+1},p^{n+1})$ are computed as
follows:
\begin{itemize}
\item[(i)] Compute the terms defined on $\Omega^{n}$ $\bigg($like $\dps{\frac{1}{\delta
    t}\int_{\Omega^n} \rho \u^n \cdot \v\,d\x}$$\bigg)$ in the system~(\ref{eq:weak-NS-disc}).
\item[(ii)] Compute  $\w^{n}=(0,0,w^n)$. We consider two options:
  \begin{itemize}
  \item 
    Scheme (M1) :  explicit treatment of the displacement of the interface
    \begin{equation}\label{eq:lapl-w-m1}
      \left\{
        \begin{array}{rcll}
          -\Delta w^n &= & 0,  & \mbox{ on } \Omega_i^n, i=1,2,\\
          w^n & =& \dps{\frac{\u^n\cdot\n_{\Sigma^n}}{n_{\Sigma^n}^{3}}}, & \mbox{ on } \Sigma^n, \\
          \dps{\frac{\partial w^n}{\partial \n}} &=& 0, & \mbox{ on } \partial \Omega,
        \end{array}
      \right.
    \end{equation}
  \item Scheme (M2) : implicit treatment of the displacement of the interface
    \begin{equation}\label{eq:lapl-w-m3}
      \left\{
        \begin{array}{rcll}
          -\Delta w^n &= & 0,  & \mbox{ on } \Omega_i^{n+1}, i=1,2,\\
          w^n \circ \invalen & =& \dps{\frac{\u^{n+1} \cdot\n_{\Sigma^{n+1}}}{n_{\Sigma^{n+1}}^{3}}}, & \mbox{ on } \Sigma^{n+1}, \\
          \dps{\frac{\partial w^n}{\partial \n}} &=& 0, & \mbox{ on } \partial \Omega.
        \end{array}
      \right.
    \end{equation}
  \end{itemize}
    where the components of the normal $\n_\Sigma$ are denoted
    by $(n_\Sigma^1,n_\Sigma^2,n_\Sigma^3)$.
  \item[(iii)] Move the nodes of the mesh according to $\alen$ defined by~(\ref{eq:alen}).
  \item[(iv)] Compute the remaining terms (defined on $\Omega^{n+1}$) in
    the system~(\ref{eq:weak-NS-disc}).
  \item[(v)] Solve (\ref{eq:weak-NS-disc}) to determine
    $(\u^{n+1},p^{n+1})$. The resolution is typically performed by a GMRES
    iterative procedure with an ILU preconditioner and
    $(\u^{n},p^{n})$ as the initial guess.
  \end{itemize}
In step (ii), the implementation of the Dirichlet boundary condition
on $w$ is made easier by defining the normals $\n_\Sigma$ at each node
of the discretized surface $\Sigma^n$. Such a definition is delicate,
since $\Sigma^n$ is only piecewise smooth, and the nodes are typically
singular points of~$\Sigma^n$. In practice,
following~\cite{engelman-sani-gresho-82}, we use approximated
normals~$\n_{\Sigma,h}$ at each node of the interface, by requiring
that the Stokes integration by parts formula holds at the discrete
level. This is one of the ingredient which ensures the exact mass
conservation of each fluids on the discretized system. We refer to
Section~5.1.3.2 in~\cite{gerbeau-le-bris-lelievre-book} or
to~\cite{gerbeau-le-bris-lelievre-03} for more details.

Note that scheme (M2) is much more expensive than scheme (M1) since it
requires the resolution of a nonlinear system (the velocity $\u^{n+1}$
and the domains $\Om^{n+1}_i$ being unknown when $\w^n$ is
computed). The resulting nonlinear system can be typically solved with
a Newton algorithm (see for example \cite{soulaimani-saad-98}). For
this study, we chose a relaxed fixed-point method, solving
iteratively systems (\ref{eq:weak-NS-disc}) and
(\ref{eq:lapl-w-m3}) at each time step.

\section{Surface Geometric Conservation Law}\label{sec:surf-gcl}

We establish in this section a few technical results that will be
useful to study the discrete energy of the system.

Let $\phi$ be a test function such that $\hat \phi$ (defined by
$\hat{\phi}(t,\hat{\x})=\phi(t,\ale(\hat{\x}))$) is
\emph{time-independent}.  Taking $\psi=1$ in~(\ref{eq:reynolds}), we
readily have
\begin{equation}
  \label{eq:gcl-cont}
  \frac{d}{dt} \int_{\Omega} \phi(t,\x)\, d\x = \int_{\Omega}
  \phi(t,\x)\, \diverg \w (t,\x) \, d\x,
\end{equation}

There is an analogous formula to compute the time derivative of the
\emph{surface} integral of $\phi$ (see~\cite{delfour-zolesio-01},
formula~(4.17) p.~355):
\begin{equation}
  \label{eq:surf-gcl-cont}
\frac{d}{dt} \int_{\Sigma_t} \phi(t,.) \,d
\sigma_{\Sigma_t}=\int_{\Sigma_t} \phi \, \tr(\nabla_{\Sigma_t} (\w)) \,d
\sigma_{\Sigma_t}.  
\end{equation}

The so-called \emph{Geometric Conservation Law} (GCL) is related to
equation~(\ref{eq:gcl-cont}). The precise definition of GCL may differ
from an author to another (see
\emph{e.g. }\cite{nkonga-guillard-94,lesoinne-farhat-95,guillard-farhat-00,formaggia-nobile-99}).
Here we will say that the GCL is satisfied when the equality in
(\ref{eq:gcl-cont}) is preserved after time discretization. The
following lemma states that the GCL is satisfied as soon as the mesh
velocity has only one nonvanishing component.
\begin{lem}[GCL]\label{lem:gcl}
  Suppose that the domain velocity $\w^n$ has the form $(0,0,w^n)$.  Let
  $\phi$ be a function defined on $\Omega_i^{n+1}$, for $i=1$ or $2$. Then
  the ALE scheme satisfies the GCL in the following sense:
    \begin{eqnarray}
\lefteqn{      \int_{\Omega_i^{n+1}}\!\!\!\!\!
      \phi(\x) \,d\x - \int_{\Omega_i^n}\!\!\! \phi\circ\alen(\y)\,
      d\y}\nonumber \\
 &=& 
\delta t \int_{\Omega_i^n} \!\!\! \phi \circ \alen(\y) \diverg_\y \w^n(\y) \, d\y, \label{eq:gcl}
\\
&=&\delta t \int_{\Omega_i^{n+1}} \!\!\!\!\!\phi (\x) \diverg_\x \left(
  \w^n\circ\invalen(\x) \right) \, d\x.      \label{eq:gcl2}
    \end{eqnarray}
\end{lem}
We refer to~\cite{gerbeau-le-bris-lelievre-03} or
to~\cite{gerbeau-le-bris-lelievre-book} for the proof of this
result. Formulas (\ref{eq:gcl}) and (\ref{eq:gcl2}) are useful to
establish the following lemma which will be used to study the effect
of the gravity on the energy balance.
\begin{lem}\label{lem:gcl_grav}
  Under the assumption of Lemma~\ref{lem:gcl}, we have:
  \begin{align}
    \int_{\Sigma^{n}} g x_3 \, \delta \rho \,  \w^{n} \cdot
    \n_{\Sigma^{n}}& =\int_{\Sigma^{n}} g x_3 \, \delta \rho \,  w^{n}
    n^3_{\Sigma^{n}}, \nonumber \\
    &= - \frac{1}{\dt}\left( \int_{\Omega^{n+1}} \rho g x_3 d\x -
      \int_{\Omega^{n}} \rho g x_3 d\x \right) - \frac{\dt}{2}
    \int_{\Sigma^{n}} \delta \rho g (w^{n})^2 n^3_{\Sigma^{n}},\label{eq:gcl1_grav}
  \end{align}
  and
  \begin{align}
    \int_{\Sigma^{n+1}} g x_3 \, \delta \rho \,&  \w^{n} \circ \invalen \cdot
    \n_{\Sigma^{n+1}} =\int_{\Sigma^{n+1}} g x_3 \, \delta \rho \,  w^{n}\circ \invalen
    n^3_{\Sigma^{n+1}}, \nonumber \\
    &= - \frac{1}{\dt}\left( \int_{\Omega^{n+1}} \rho g x_3 d\x -
      \int_{\Omega^{n}} \rho g x_3 d\x \right) + \frac{\dt}{2}
    \int_{\Sigma^{n}} \delta \rho g (w^{n})^2 n^3_{\Sigma^{n}}.\label{eq:gcl2_grav}
  \end{align}
\end{lem}
\begin{proof}
  Formula (\ref{eq:gcl}) applied to the function $\rho g x_3$ yields
\begin{align*}
\int_{\Omega^{n+1}} \rho g x_3 
&=\int_{\Omega^{n}} \rho g (y_3 + \dt w^n(\y)) \, d\y + \dt
\int_{\Omega^{n}} \rho g (y_3 + \dt w^n(\y)) \partial_{y_3} w^n \, d\y,\\
&=\int_{\Omega^{n}} \rho g y_3 + \dt \left( \int_{\Omega^{n}} \rho
  gw^n(\y) \, d\y + \int_{\Omega^{n}} \rho g y_3 \partial_{y_3} w^n \, d\y
\right) + \dt^2 \int_{\Omega^{n}} \rho g w^n \partial_{y_3} w^n \, d\y,\\
&=\int_{\Omega^{n}} \rho g y_3 + \dt \int_{\Omega^{n}} \rho g
\partial_{y_3} (y_3 w^n) \, d\y + \frac{\dt^2}{2} \int_{\Omega^{n}} \rho g \partial_{y_3} ((w^n)^2) \, d\y,\\
&=\int_{\Omega^{n}} \rho g y_3 - \dt \int_{\Sigma^{n}} \delta \rho g
y_3 w^n n^3_{\Sigma^{n}}  - \frac{\dt^2}{2} \int_{\Sigma^{n}} \delta
\rho g (w^n)^2 n^3_{\Sigma^{n}},
\end{align*}
from which we deduce (\ref{eq:gcl1_grav}). We recall that $\delta
\rho$ is defined by~\eqref{eq:delta_rho}. For the second formula, we
use~(\ref{eq:gcl2}):
\begin{align*}
\int_{\Omega^{n+1}} \rho g x_3 
&=\int_{\Omega^{n}} \rho g (y_3 + \dt w^n(\y)) \, d\y + \dt
\int_{\Omega^{n+1}} \rho g x_3 \partial_{x_3} \left( w^n \circ \invalen
  (\x) \right)\, d\x,\\
&=\int_{\Omega^{n}} \rho g (y_3 + \dt w^n(\y)) \, d\y - \dt
\int_{\Omega^{n+1}} \partial_{x_3} \left( \rho g x_3 \right)  w^n \circ \invalen
  (\x) \, d\x,\\
&=\int_{\Omega^{n}} \rho g (y_3 + \dt w^n(\y)) \, d\y - \dt
\int_{\Omega^{n+1}}  \rho g   w^n \circ \invalen
  (\x) \, d\x  \\
&\quad - \dt
\int_{\Sigma^{n+1}} \delta \rho  g x_3 n^3_{\Sigma^{n+1}}  w^n \circ \invalen,\\
&=\int_{\Omega^{n}} \rho g y_3 + \dt \left( \int_{\Omega^{n}} \rho gw^n(\y)) \, d\y -
\int_{\Omega^{n+1}}  \rho g   w^n \circ \invalen
  (\x) \, d\x \right)  \\
&\quad - \dt
\int_{\Sigma^{n+1}} \delta \rho  g x_3 n^3_{\Sigma^{n+1}}  w^n \circ \invalen,\\
&=\int_{\Omega^{n}} \rho g y_3 - \dt^2 \left( \int_{\Omega^{n}} \rho g
  w^n\partial_{y_3} w^n \, d\y \right)  - \dt
\int_{\Sigma^{n+1}} \delta \rho  g x_3 n^3_{\Sigma^{n+1}}  w^n \circ \invalen,\\
&=\int_{\Omega^{n}} \rho g y_3 - \dt
\int_{\Sigma^{n+1}} \delta \rho  g x_3 n^3_{\Sigma^{n+1}}  w^n \circ \invalen - \frac{\dt^2}{2} \left( \int_{\Omega^{n}} \rho g
  \partial_{y_3} (w^n)^2  \right), \\
&=\int_{\Omega^{n}} \rho g y_3 - \dt
\int_{\Sigma^{n+1}} \delta \rho  g x_3 n^3_{\Sigma^{n+1}}  w^n \circ \invalen + \frac{\dt^2}{2} \left( \int_{\Sigma^{n}} \delta\rho \,  g
   (w^n)^2 n^3_{\Sigma^{n}}  \right),
\end{align*}
which gives (\ref{eq:gcl2_grav}).
\end{proof}

The above results are useful to study the energy balance at the discrete
level in the presence of gravity. To study the effect of surface
tension on the energy balance,
we have to introduce a notion of "\emph{Surface GCL}", namely we have
to check how equality (\ref{eq:surf-gcl-cont}) is preserved by the
time scheme. Contrarily to the volume case, we are not able to prove
that the scheme satisfies a "Surface GCL". Nevertheless, we can
establish inequalities which will be convenient in the sequel.
\begin{lem}[Surface GCL]\label{lem:gcl_surf}
Suppose that the domain velocity $\w^n$ has the form $(0,0,w^n)$. Let
$\phi$ be a function defined on $\Sigma^{n+1}$. Then, if $\dt$ is
sufficiently small so that
\begin{equation}\label{eq:hyp_dt_small_1}
1 + \dt \, \tr \left( \nabla_{\Sigma^{n}}
\left(\w^{n} \right) \right) \geq 0 \qquad \text{on $\Sigma^n$},
\end{equation}
the ALE scheme is such that:
\begin{equation}
\int_{\Sigma^{n+1}} \phi \, d \sigma_{\Sigma^{n+1}} - \int_{\Sigma^{n}}
\phi \circ \alen \, d \sigma_{\Sigma^{n}} \geq \dt \int_{\Sigma^n} \phi
\circ \alen \tr(\nabla_{\Sigma^n}(\w^n))\, d
\sigma_{\Sigma^{n}}.\label{eq:gcl_surf_1}
\end{equation}
Likewise,  if $\dt$ is
sufficiently small so that
\begin{equation}\label{eq:hyp_dt_small_2}
1 - \dt \, \tr \left( \nabla_{\Sigma^{n+1}}
\left(\w^{n}  \circ
  \invalen \right) \right) \geq 0\qquad \text{on $\Sigma^{n+1}$},
\end{equation}
the ALE scheme is such that:
\begin{equation}
\int_{\Sigma^{n+1}} \phi\, d \sigma_{\Sigma^{n+1}} - \int_{\Sigma^{n}}
\phi \circ \alen\, d \sigma_{\Sigma^{n}} \leq \dt \int_{\Sigma^{n+1}} \phi
\, \tr(\nabla_{\Sigma^{n+1}} (\w^n \circ \invalen))\, d
\sigma_{\Sigma^{n+1}}.\label{eq:gcl_surf_2}
\end{equation}
Moreover, in both cases, the difference between the two sides of the
inequalities~(\ref{eq:gcl_surf_1}) and (\ref{eq:gcl_surf_2}) is of
order $\dt^2$ in the limit $\dt \to 0$.
\end{lem}

\begin{proof}
We only consider~(\ref{eq:gcl_surf_1}), the proof
of~(\ref{eq:gcl_surf_2}) being similar. By a change of variable we have
(see~\cite[Formula~(4.9) p.~353]{delfour-zolesio-01}):
\begin{align*}
\int_{\Sigma^{n+1}} \phi d \sigma_{\Sigma^{n+1}}
&=\int_{\Sigma^{n}} \phi \circ \alen \left| {\rm Cof}(\nabla \alen)
  \n_{\Sigma^{n}} \right| d \sigma_{\Sigma^{n}},\\
&=\int_{\Sigma^{n}}  \phi \circ \alen \left|\n_{\Sigma^{n}} + \dt \left(\begin{array}{c}
  n^1_{\Sigma^{n}} \partial_{y_3} w^n  -
  n^3_{\Sigma^{n}} \partial_{y_1} w^n  \\ n^2_{\Sigma^{n}} \partial_{y_3} w^n  -
  n^3_{\Sigma^{n}} \partial_{y_2} w^n  \\ 0 \end{array}  \right)  \right| d \sigma_{\Sigma^{n}},
\end{align*}
where ${\rm Cof}(\nabla \alen)$ denotes the matrix of cofactors of the
Jacobian $\nabla \alen$. 

Let us consider the difference
\begin{equation}\label{eq:3}
\left|\n_{\Sigma^{n}} + \dt \left(\begin{array}{c}
  n^1_{\Sigma^{n}} \partial_{y_3} w^n   -
  n^3_{\Sigma^{n}} \partial_{y_1} w^n   \\ n^2_{\Sigma^{n}} \partial_{y_3} w^n  -
  n^3_{\Sigma^{n}} \partial_{y_2} w^n  \\ 0 \end{array}  \right)
\right|^2 - \left| 1 + \dt \, \tr \left( \nabla_{\Sigma^{n}}
\left(\w^{n} \right) \right) \right|^2.
\end{equation}
The first term writes:
\begin{align*}
&1 + 2 \dt \left( (n^1_{\Sigma^{n}})^2 \partial_{y_3} w^n  -
  n^1_{\Sigma^{n}} n^3_{\Sigma^{n}} \partial_{y_1} w^n +  (n^2_{\Sigma^{n}})^2 \partial_{y_3} w^n  -
  n^2_{\Sigma^{n}}n^3_{\Sigma^{n}} \partial_{y_2} w^n  \right)\\
&+ \dt^2 \left( \left( n^1_{\Sigma^{n}} \partial_{y_3} w^n  -
  n^3_{\Sigma^{n}} \partial_{y_1} w^n \right)^2 + \left(
  n^2_{\Sigma^{n}} \partial_{y_3} w^n  - n^3_{\Sigma^{n}} \partial_{y_2} w^n  \right)^2 \right).
\end{align*}
For the second term, we have:
\begin{align*}
\tr \left( \nabla_{\Sigma^{n}}
\left(\w^{n}\right) \right)
&=\diverg(\w^n) - \nabla \w^n \n_{\Sigma^{n}} \cdot \n_{\Sigma^{n}}, \\
&= \partial_{y_3} w^n - n^3_{\Sigma^{n}} \left( \sum_{i=1}^3 n^i_{\Sigma^{n}} \partial_{y_i}
    w^n\right),\\
&= n^1_{\Sigma^{n}}\left( n^1_{\Sigma^{n}}\partial_{y_3}
    w^n - n^3_{\Sigma^{n}}\partial_{y_1}
    w^n \right)+n^2_{\Sigma^{n}}\left( n^2_{\Sigma^{n}}\partial_{y_3}
    w^n - n^3_{\Sigma^{n}}\partial_{y_2}
    w^n \right).
\end{align*}
The second term in~(\ref{eq:3}) thus writes:
\begin{align*}
&1 + 2 \dt \left(n^1_{\Sigma^{n}}\left( n^1_{\Sigma^{n}}\partial_{y_3}
    w^n - n^3_{\Sigma^{n}}\partial_{y_1}
    w^n \right)+n^2_{\Sigma^{n}}\left( n^2_{\Sigma^{n}}\partial_{y_3}
    w^n - n^3_{\Sigma^{n}}\partial_{y_2}
    w^n \right) \right)\\
& + \dt^2 \left( n^1_{\Sigma^{n}}\left( n^1_{\Sigma^{n}}\partial_{y_3}
    w^n - n^3_{\Sigma^{n}}\partial_{y_1}
    w^n \right)+n^2_{\Sigma^{n}}\left( n^2_{\Sigma^{n}}\partial_{y_3}
    w^n - n^3_{\Sigma^{n}}\partial_{y_2}
    w^n \right) \right)^2.
\end{align*}
Using now Cauchy-Schwarz inequality, we obtain:
\begin{align}
&\left( n^1_{\Sigma^{n}}\left( n^1_{\Sigma^{n}}\partial_{y_3}
    w^n - n^3_{\Sigma^{n}}\partial_{y_1}
    w^n \right)+n^2_{\Sigma^{n}}\left( n^2_{\Sigma^{n}}\partial_{y_3}
    w^n - n^3_{\Sigma^{n}}\partial_{y_2}
    w^n \right) \right)^2 \nonumber \\
&\leq \left( (n^1_{\Sigma^{n}})^2 + (n^2_{\Sigma^{n}})^2 \right)
\left( \left( n^1_{\Sigma^{n}} \partial_{y_3} w^n  -
  n^3_{\Sigma^{n}} \partial_{y_1} w^n \right)^2 + \left(
  n^2_{\Sigma^{n}} \partial_{y_3} w^n  - n^3_{\Sigma^{n}}
  \partial_{y_2} w^n  \right)^2 \right), \nonumber \\
&\leq \left( n^1_{\Sigma^{n}} \partial_{y_3} w^n  -
  n^3_{\Sigma^{n}} \partial_{y_1} w^n \right)^2 + \left(
  n^2_{\Sigma^{n}} \partial_{y_3} w^n  - n^3_{\Sigma^{n}}
  \partial_{y_2} w^n  \right)^2.\label{eq:4}
\end{align}
This shows that~(\ref{eq:3}) is nonnegative and bounded from above by
$C \dt^2 \left\| \frac{\partial \w^n}{\partial \y} \right\|^2$,
where~$C$ denotes a constant.

Thus, if $\dt$ is sufficiently small so that~(\ref{eq:hyp_dt_small_1})
is satisfied, we have:
\begin{align*}
\int_{\Sigma^{n+1}} \phi d \sigma_{\Sigma^{n+1}}
&\geq \int_{\Sigma^{n}} \phi \circ \alen \left( 1 + \dt \, \tr \left( \nabla_{\Sigma^{n}}
\left(\w^{n} \right) \right) \right) d \sigma_{\Sigma^{n}},\\
&\geq  \int_{\Sigma^{n}} \phi \circ \alen d \sigma_{\Sigma^{n}} + \dt
\int_{\Sigma^{n}}  \phi \circ \alen \tr \left( \nabla_{\Sigma^{n}}
\left(\w^{n} \right) \right)d \sigma_{\Sigma^{n}}.
\end{align*}

Moreover, if $\dt$ is
sufficiently small so that $1 + \dt \, \tr \left( \nabla_{\Sigma^{n}}
\left(\w^{n} \right) \right) \geq \epsilon > 0$ on $\Sigma^n$, then
the  difference between the two sides of the
inequality~(\ref{eq:gcl_surf_1}) writes:
\begin{align}
0 &\leq
\int_{\Sigma^{n+1}} \phi \, d \sigma_{\Sigma^{n+1}} - \int_{\Sigma^{n}}
\phi \circ \alen \, d \sigma_{\Sigma^{n}} - \dt \int_{\Sigma^n} \phi
\circ \alen \tr(\nabla_{\Sigma^n}(\w^n))\, d
\sigma_{\Sigma^{n}},\nonumber \\
&= \int_{\Sigma^n} \phi \circ \alen \left( \left|\n_{\Sigma^{n}} + \dt \left(\begin{array}{c}
  n^1_{\Sigma^{n}} \partial_{y_3} w^n  -
  n^3_{\Sigma^{n}} \partial_{y_1} w^n  \\ n^2_{\Sigma^{n}} \partial_{y_3} w^n  -
  n^3_{\Sigma^{n}} \partial_{y_2} w^n  \\ 0 \end{array}  \right)
\right| - |1 + \dt \, \tr(\nabla_{\Sigma^n}(\w^n))|  \right) d
\sigma_{\Sigma^{n}}. \nonumber
\end{align}
Using now the inequality $0 \le |\uu{a}|-|b| \le \frac{|\uu{a}|^2-|b|^2}{2|b|}$ with $\uu{a}=\n_{\Sigma^{n}} + \dt \left(\begin{array}{c}
  n^1_{\Sigma^{n}} \partial_{y_3} w^n  -
  n^3_{\Sigma^{n}} \partial_{y_1} w^n  \\ n^2_{\Sigma^{n}} \partial_{y_3} w^n  -
  n^3_{\Sigma^{n}} \partial_{y_2} w^n  \\ 0 \end{array}  \right)$ and
$b=1 + \dt \, \tr(\nabla_{\Sigma^n}(\w^n))$, we obtain the following
estimate of the  difference between the two sides of the
inequality~(\ref{eq:gcl_surf_1}):
\begin{align}
\int_{\Sigma^{n+1}} &\phi \, d \sigma_{\Sigma^{n+1}} - \int_{\Sigma^{n}}
\phi \circ \alen \, d \sigma_{\Sigma^{n}} - \dt \int_{\Sigma^n} \phi
\circ \alen \tr(\nabla_{\Sigma^n}(\w^n))\, d
\sigma_{\Sigma^{n}},\nonumber\\
&\leq
\Frac{C}{\epsilon} \left\| \frac{\partial \w^n}{\partial \y} \right\|^2_{L^2(\Sigma^n)} \left\|
\phi \right\|_{L^\infty(\Sigma^{n+1})}  \dt^2.\label{eq:tension-dt-2}
\end{align}

This concludes the proof.
\end{proof}
\begin{remark} By analogy with the GCL, we would say that the scheme
  satisfies the "Surface GCL" if the inequalities in
  (\ref{eq:gcl_surf_1}) or (\ref{eq:gcl_surf_2}) were
  equalities. Looking at the proof of Lemma~\ref{lem:gcl_surf}
  (see~(\ref{eq:4})), one can see that the equality
  in~(\ref{eq:gcl_surf_1}) would be obtained if and only
  if
\begin{align*}
n^3_{\Sigma^{n}}=0,
\end{align*}
which is not satisfied in practice. Devising a scheme that allows to
recover an equality in (\ref{eq:gcl_surf_1}) or (\ref{eq:gcl_surf_2})
(and keeping equalities in (\ref{eq:gcl}) and (\ref{eq:gcl2}))
is, to our knowledge, an interesting open question.
\end{remark}

\section{Energy estimates and time-discretization}\label{sec:energy}

We now investigate the stability in the energy norm of the implicit
and explicit numerical schemes. In Section~\ref{sec:energy_cont}, we
first recall the derivation of the energy estimate at the continuous
level. We then discuss in Section~\ref{sec:energy_disc} to what extent
the computation at the continuous level can be reproduced on the
time-discretized system. In particular, we will emphasize the effect
of the explicit treatment of the displacement of the free interface on the energy balance.
We introduce the following notations: 
$$
K(t) = \frac{1}{2}\int_{\Omega} \rho(t,\x) |\u(t,\x)|^2~~\mbox{(kinetic energy)},
$$
$$
P_v(t) = \int_{\Omega}
    \frac{\eta(t,\x)}{2} \left| \nabla \u(t,\x) + (\nabla \u(t,\x))^T \right|^2~~\mbox{(viscous power)},
$$
$$
W(t) = \int_{\Omega} \rho(t,\x) g x_3~~\mbox{(potential energy)},
$$
$$
\|\u\|_{\partial\Om} = \sqrt{\int_{\partial \Omega} |\u|^2}
\mbox{ and } (\u,\v)_{\partial\Om}= \int_{\partial\Om} \u\cdot\v.
$$
We denote by $|\Sigma(t)|$ the measure of the surface $\Sigma(t)$,
namely $\int_{\Sigma} d \sigma_\Sigma$.

\subsection{Energy estimate at the continous level}\label{sec:energy_cont}


\begin{prop}
The physical system described in
  Section~\ref{sec:GNBC} satisfies the following energy equality:
  \begin{equation}\label{eq:NRJ}
    \frac{dK}{dt} + \frac{dW}{dt}  + 
    \gamma
    \frac{d|\Sigma(t)|}{dt}  +  P_v
  + \beta (\u,\u-\u_b)_{\partial\Om} = \gamma \int_{\partial\Sigma} \cos(\theta_s) \t_{\partial\Om}\cdot \u\,dl_{\partial\Sigma}.
  \end{equation}
\end{prop}

\proof


The proof of this proposition is standard, we briefly sketch it for
convenience.  Multiplying the strong form of the equations by $\u$ and
integrating, we obtain:
 \begin{equation}
\label{eq:eq42}
\begin{array}{l}
  \dps{\frac{1}{2} \frac{d}{dt}   \int_\Omega \rho |\u|^2 + \int_\Omega
    \frac{\eta}{2} \left| \nabla \u + \nabla \u^T \right|^2  + 
\beta(\u-\u_b,\u)_{\partial\Om}}\\
\quad  \dps{=   -\gamma \int_{\Sigma} \tr(\nabla_\Sigma
  \u) \, d\sigma_{\Sigma} +   \int_\Omega  \rho \g \cdot \u +  \gamma \int_{\partial\Sigma} \cos(\theta_s) \t_{\partial\Om}\cdot \u\,dl_{\partial\Sigma}.}
  \end{array}
 \end{equation}
Here are the main steps to obtain this equation:
 \begin{eqnarray}
   \label{3eq:non-2}
   \int_\Omega \frac{\partial (\rho\u)}{\partial t} \cdot\u\,d\x &=&
   \int_\Omega \frac{\partial \rho}{\partial t} |\u|^2
   \,d\x+\int_\Omega\rho\,\u\cdot \frac{\partial \u}{\partial t} \,d\x,\nonumber\\
 &=&{1\over
   2}\int_\Omega\frac{\partial \rho}{\partial t}\, |\u|^2 \,d\x+{1\over
   2}{d\over{dt}}\int_\Omega\rho |\u|^2\,d\x\nonumber,\\
 &=&{1\over 2}\int_\Omega\rho\,\u\cdot\nabla (\u^2)\,d\x+{1\over
   2}{d\over{dt}}\int_\Omega\rho|\u|^2\,d\x. \nonumber
 \end{eqnarray}
 The first term in the right-hand side of the last equality cancels
 with the term resulting from the integration of $\diverg(\rho
 \u\otimes\u) \cdot \u$. The remaining terms of (\ref{eq:eq42}) comes from the
 integration by parts of the stress terms and the GNBC (\ref{eq:GNBC}).

 We now consider the surface tension term. We have, since $\w \cdot
 \n_{\partial \Omega} =0$ (see~\cite[Formula~(4.17) p.~355]{delfour-zolesio-01}),
 \begin{align}
\int_{\Sigma} \tr (\nabla_\Sigma \u) d \sigma_\Sigma
= \int_{\Sigma} \tr(\nabla_\Sigma \w) d \sigma_\Sigma
= \frac{d}{dt} \int_{\Sigma} d \sigma_\Sigma.  \label{eq:NRJ_TS}
 \end{align}
For the last equality, we refer for example
to~\cite{delfour-zolesio-01}, formula (4.17) p. 355. The first equality relies
on~(\ref{eq:div_th}) and on the fact that $\u \cdot \n_\Sigma=\w \cdot
\n_\Sigma$ which implies that $\u \cdot \m=\w \cdot \m$.

Concerning the gravity term, we have:
\begin{align}
\int_\Omega \rho \g \cdot \u
&= - \int_\Omega \rho g \nabla x_3 \cdot \u = \int_\Omega \diverg(\rho \u) g x_3, \nonumber \\
&= - \int_\Omega \frac{\partial \rho}{\partial t} g x_3 = - \frac{d}{d t} \int_\Omega \rho g x_3, \label{eq:NRJ_grav}
\end{align}
where we have used the fact that 
$\nabla \rho = \delta \rho \, \n_{\Sigma} \, \delta_{\Sigma}$
(see~(\ref{eq:delta_rho}) for the
definition of $\delta \rho$).
\endproof

\subsection{Energy estimates at the discrete level}\label{sec:energy_disc}

In this section, we present the discrete counterpart of the
computations previously made for the continuous system. For
simplicity, the problem is only discretized in time.  Nevertheless,
all the following computations could be carried out after space
discretization with finite elements as soon as the pressure finite
element space contains linear functions and with a definition of the
normals at the interface consistent with the discrete Stokes formula (see~\cite{gerbeau-le-bris-lelievre-03,gerbeau-le-bris-lelievre-book}).

We define the discrete counterpart of the notations introduced before:
$$
K^n = \frac{1}{2}\int_{\Omega^{n}} \rho |\u^{n}|^2~~\mbox{(kinetic energy)},
$$
$$
P_v^n = \int_{\Omega^{n}}
\frac{\eta}{2} \left| \nabla \u^{n} + (\nabla \u^{n})^T \right|^2~~\mbox{(viscous power)},
$$
$$
W^n = \int_{\Omega^{n}} \rho g x_3~~\mbox{(potential energy)},
$$
and we denote by $|\Sigma^n|$ the measure of the surface
$\Sigma^n$.

\begin{prop}
\label{prop:balance-1}
  We first consider the case without body force and without surface tension,
  namely we assume $g=0$ and $\gamma=0$. We suppose that the domain
  movement is computed with the explicit scheme (M1), \emph{i.e.} solving
  equation (\ref{eq:lapl-w-m1}). Then, we have the following stability result:
  \begin{equation}
    \frac{K^{n+1}-K^n}{\dt} + P_v^{n+1}  + \beta(\u^{n+1}-\u_b,\u^{n+1})_{\partial\Om} =
    - \int_{\Omega^n}  \frac{\rho}{2 \dt} |\u^{n+1} \circ \alen - \u^n|^2.
\label{eq:NRJ_disc}  \end{equation}  
\end{prop}
The estimate~(\ref{eq:NRJ_disc}) is the discrete counterpart
of~(\ref{eq:NRJ}).  For the proof, we refer the reader to
\cite{gerbeau-le-bris-lelievre-03} where a similar result is
presented.  It is based in particular on the GCL (\ref{eq:gcl}).
Since the right-hand side is nonpositive, it shows that the time
discretization scheme does not bring spurious energy to the system.
Such a property is easy to obtain on a fixed mesh, but more
complicated to ensure on a moving mesh
(see~\cite{gerbeau-le-bris-lelievre-03}
or~\cite{gerbeau-le-bris-lelievre-book} for the details).  It is
noteworthy that such a stability property is obtained with a scheme
which only requires to solve a linear system at each time step (since
scheme (M1) decouples the mesh movement and the fluid
resolution). Nevertheless, we have to bear in mind that the assumption
$g=0$ and $\gamma=0$ is not satisfied in the interesting physical
configurations. We will show in the sequel that when gravity or
surface tension are considered, an analogous stability in the energy
norm can be obtained with an \emph{implicit} treatment of the domain
movement (scheme (M2)) but is no more valid with an explicit scheme
like (M1).

\begin{prop}
\label{prop:balance-2}
We suppose that the domain movement is computed with the explicit
scheme (M1), \emph{i.e.} solving equation (\ref{eq:lapl-w-m1}) and
that the gravity is computed on $\Om^{n+1}$ (namely $\Om^\star=
\Om^{n+1}$ in \eqref{eq:weak-NS-disc}). Then we have the following
stability result:
  \begin{equation}
    \begin{array}[t]{l}
      \Frac{K^{n+1} - K^n}{\delta t} + 
      \Frac{W^{n+2} - W^{n+1}}{\delta t} + P_v^{n+1}
+ \gamma\Frac{|\Sigma^{n+2}| - |\Sigma^{n+1}|}{\dt} 
+ \beta(\u^{n+1}-\u_b,\u^{n+1})_{\partial\Om}
      \\[9pt]   ~~~~~~~~+ \Int{\Om^n}{}\Frac{\rho}{2\delta t}|\u^{n+1} - \u^n|^2 
      = 
\epsilon_g^{n+1} + \epsilon_{\gamma,exp}^{n+1} + \gamma \int_{\partial\Sigma^{n+1}} \cos(\theta_s) \t_{\partial\Om}\cdot \u^{n+1}\,dl_{\partial\Sigma}
    \end{array}
    \label{eq:NRJ_disc_2}
  \end{equation}  
with
\begin{equation}
  \label{eq:epsilon-g}
\epsilon_{g}^{n+1} = - \Frac{\delta t}{2}\Int{\Sigma^{n+1}}{} \delta \rho g
      (w^{n+1})^2 n_{\Sigma^{n+1},3}  \geq 0,
\end{equation}
and
\begin{equation}
  \label{eq:epsilon-gamma}
\epsilon_{\gamma,exp}^{n+1} = \frac{\gamma}{\dt} \left( |\Sigma^{n+2}| - |\Sigma^{n+1}|- \dt \int_{\Sigma^{n+1}} \tr(\nabla_{\Sigma^{n+1}}
  \w^{n+1}) \right) \geq 0.
\end{equation}
\end{prop}

\begin{remark}
  In the above Proposition, note that the terms $\epsilon_g^n$ and
  $\epsilon_{\gamma,exp}^n$ are of order $\delta t$, but are
  nonnegative: they bring a spurious power to the system. This is a
  consequence of the explicit treatment of the free interface displacement.
\end{remark}
\proof 

We only focus on the terms related to gravity and surface tension. The
other terms can be treated as in Proposition~\ref{prop:balance-1},
following the arguments of \cite{gerbeau-le-bris-lelievre-03}.

Consider first the gravity term. The purpose is to mimic at the
discrete level the computations done to
establish~(\ref{eq:NRJ_grav}). Multiplying the equation by $\u^{n+1}$,
we have:
\begin{align*}
\int_{\Omega^{n+1}} \rho \g \cdot \u^{n+1}
&=-\int_{\Omega^{n+1}} \rho g \nabla x_3 \cdot  \u^{n+1},\\
&=\int_{\Omega^{n+1}}  g x_3  \diverg( \rho  \u^{n+1}), \\
&=\int_{\Sigma^{n+1}} g x_3 \, \delta \rho \,  \u^{n+1} \cdot \n_{\Sigma^{n+1}}.
\end{align*}
Therefore, noticing that with Scheme (M1) we have $\u^{n+1} \cdot
\n_{\Sigma^{n+1}} = \w^{n+1}\cdot \n_{\Sigma^{n+1}}$, and using
formula~\eqref{eq:gcl1_grav} of Lemma~\ref{lem:gcl_grav} we have:
\begin{align}
\int_{\Omega^{n+1}} \rho \g \cdot \u^{n+1}
&=\int_{\Sigma^{n+1}} g x_3 \, \delta \rho \,  \w^{n+1} \cdot
\n_{\Sigma^{n+1}}, \nonumber \\
&=\int_{\Sigma^{n+1}} g x_3 \, \delta \rho \,  w^{n+1}
n^3_{\Sigma^{n+1}}, \nonumber \\
&= - \frac{1}{\dt}\left( \int_{\Omega^{n+2}} \rho g x_3 d\x -
  \int_{\Omega^{n+1}} \rho g x_3 d\x \right) \nonumber \\
& \quad - \frac{\dt}{2}
\int_{\Sigma^{n+1}} \delta \rho g (w^{n+1})^2 n^3_{\Sigma^{n+1}}.\label{eq:NRJ_grav_disc}
\end{align}
We note that the last term in~(\ref{eq:NRJ_grav_disc}), namely $\epsilon_g^{n+1}$, is of
order~$\delta t$ and it may bring some energy in the system
since it is nonnegative. Indeed, if the heaviest fluid is below
the lightest one (which is the case in our practical applications),
$\delta \rho \leq 0$ and $n^3_{\Sigma^{n+1}} \geq 0$.

We next consider the surface tension term and we try to reproduce on
the discrete system the formula~(\ref{eq:NRJ_TS}). Using again $\u^{n+1}
\cdot \n_{\Sigma^{n+1}} = \w^{n+1}\cdot \n_{\Sigma^{n+1}}$, we have with~(\ref{eq:div_th}),
\begin{align}
- \gamma \int_{\Sigma^{n+1}} \tr(\nabla_{\Sigma^{n+1}} \u^{n+1})
&=- \gamma \int_{\Sigma^{n+1}} \tr(\nabla_{\Sigma^{n+1}}
\w^{n+1}), \nonumber \\
&=- \frac{\gamma}{\dt} \left(|\Sigma^{n+2}| - |\Sigma^{n+1}|\right) \nonumber \\
& \quad + \frac{\gamma}{\dt} \left( |\Sigma^{n+2}| - |\Sigma^{n+1}|- \dt \int_{\Sigma^{n+1}} \tr(\nabla_{\Sigma^{n+1}}
  \w^{n+1}) \right) \label{eq:NRJ_TS_disc}.
\end{align}
To estimate the last term in~(\ref{eq:NRJ_TS_disc}), namely $\epsilon_{\gamma,exp}^{n+1} $, we need the Surface GCL result proved in
Lemma~\ref{lem:gcl_surf}.  By choosing $\phi=1$
in~(\ref{eq:gcl_surf_1}), we see that the last term
in~(\ref{eq:NRJ_TS_disc}) is of order $\delta t$ and that it is
nonnegative in the limit $\dt \to 0$ (see (\ref{eq:tension-dt-2})).

\endproof

\begin{prop}
  \label{prop:balance-3}
  We suppose that the domain movement is computed with the implicit
  scheme (M2), \emph{i.e.} solving equation (\ref{eq:lapl-w-m1}), and
  that the gravity is computed on $\Om^{n+1}$ (namely $\Om^\star=
  \Om^{n+1}$ in \eqref{eq:weak-NS-disc}). Then we have the following
  stability result:
  \begin{equation}
    \begin{array}[t]{l}
      \Frac{K^{n+1} - K^n}{\delta t} + 
      \Frac{W^{n+1} - W^{n}}{\delta t} + P_v^{n+1}
+ \gamma\Frac{|\Sigma^{n+1}| - |\Sigma^{n}|}{\dt}
+ \beta(\u^{n+1}-\u_b,\u^{n+1})_{\partial\Om} 
      \\[9pt]   ~~~~~~~~+ \Int{\Om^n}{}\Frac{\rho}{2\delta t}|\u^{n+1} - \u^n|^2 
      = 
 - \epsilon_g^n  - \epsilon_{\gamma,imp}^n  + \gamma \int_{\partial\Sigma^{n+1}} \cos(\theta_s) \t_{\partial\Om}\cdot \u^{n+1}\,dl_{\partial\Sigma}
    \end{array}
    \label{eq:NRJ_disc_2b}
  \end{equation}  
with
\begin{equation}
  \label{eq:epsilon-g-impl}
  \epsilon_g^n = - \Frac{\delta t}{2}\Int{\Sigma^n}{} \delta \rho g (w^n)^2 n_{\Sigma^n,3} \geq 0,
\end{equation}
and
\begin{equation}
  \label{eq:epsilon-gamma-impl}
\epsilon_{\gamma,imp}^{n} = 
- \frac{\gamma}{\dt} \left( |\Sigma^{n+1}| - |\Sigma^{n}|- \dt \int_{\Sigma^{n+1}} \tr(\nabla_{\Sigma^{n+1}} \w^{n}) \right) \geq 0.
\end{equation}
The scheme is therefore stable in the energy norm.
\end{prop}
\proof

Consider first the gravity term. With the implicit scheme (M2), we
have $\w^n\cdot\n_{\Sigma^{n+1}} = \u^{n+1} \cdot\n_{\Sigma^{n+1}}$
(see (\ref{eq:lapl-w-m3})). Thus, using formula~\eqref{eq:gcl2_grav}
of Lemma~\ref{lem:gcl_grav}, the counterpart of~(\ref{eq:NRJ_grav})
reads (compare with~\eqref{eq:NRJ_grav_disc}):
\begin{align}
\int_{\Omega^{n+1}} \rho \g \cdot \u^{n+1}
&=\int_{\Sigma^{n+1}} g x_3 \, \delta \rho \,  \u^{n+1} \cdot
\n_{\Sigma^{n+1}}, \nonumber \\
&=\int_{\Sigma^{n+1}} g x_3 \, \delta \rho \,  w^{n}\circ \invalen
n^3_{\Sigma^{n+1}}, \nonumber \\
&= - \frac{1}{\dt}\left( \int_{\Omega^{n+1}} \rho g x_3 d\x -
  \int_{\Omega^{n}} \rho g x_3 d\x \right) \nonumber \\
& \quad + \frac{\dt}{2}
\int_{\Sigma^{n}} \delta \rho g (w^{n})^2 n^3_{\Sigma^{n}}.\label{eq:NRJ_grav_disc_impl}
\end{align}
The last term is exactly $\epsilon_g^n$. Consider now the surface tension term:
\begin{align}
  - \gamma \int_{\Sigma^{n+1}} \tr(\nabla_{\Sigma^{n+1}} \u^{n+1})
  &=- \gamma \int_{\Sigma^{n+1}} \tr(\nabla_{\Sigma^{n+1}}
  \w^{n}), \nonumber \\
  &=- \frac{\gamma}{\dt} \left(|\Sigma^{n+1}| - |\Sigma^{n}|\right) \nonumber \\
  & \quad + \frac{\gamma}{\dt} \left( |\Sigma^{n+1}| - |\Sigma^{n}|- \dt \int_{\Sigma^{n+1}} \tr(\nabla_{\Sigma^{n+1}}
    \w^{n}) \right) \label{eq:NRJ_TS_disc_impl}.
\end{align}
To estimate the last term in~(\ref{eq:NRJ_TS_disc}), namely $-
\epsilon_{\gamma,imp}^{n} $, we use again Lemma~\ref{lem:gcl_surf}
(surface GCL).  By choosing $\phi=1$ in~(\ref{eq:gcl_surf_2}), we see
that it is nonpositive in the limit $\dt \to 0$.

\endproof

\subsection{Discussion}

In summary, we have shown that the numerical scheme with an explicit
treatment of the free surface displacement is stable in the energy norm in absence
of gravity and surface tension. But in presence of gravity and/or
surface tension, some terms may introduce a spurious energy in the
system if the treatment of the interface displacement is explicit (and if the
timestep is not sufficiently small). Moreover, we observe that in both
cases (gravity and surface tension) the terms of order $\delta t$
((\ref{eq:NRJ_grav_disc}) or (\ref{eq:NRJ_TS_disc})) which have the
wrong sign are integrals over the interface between the two fluids.
These theoretical results are in agreement with typical
observations. We indeed often notice that when a numerical instability
occurs, it is located on the interface between the two
fluids. Moreover, such instabilities typically increase with
increasing gravity or surface tension.

We have also shown that an implicit treatment of the free surface
displacement yields a numerical scheme which is stable in the energy
norm. Of course, such a scheme is more expensive, and we show by
numerical experiments in the next section how to build schemes with are
cheap but seems to have similar stability properties.


\section{Numerical experiments}\label{sec:num}

In all the following numerical experiments, the Navier-Stokes
equations are discretized with the Q2/P1 pair of finite element, with
a discontinuous pressure.

\subsection{Energy balance in presence of gravity}\label{sec:num-energy}

The purpose of the simulations presented in this section is to
illustrate on a simple physical experiment the theoretical results
established in Section~\ref{sec:energy}.
\begin{figure}[htbp]
  \centering
  \centerline{\epsfig{file=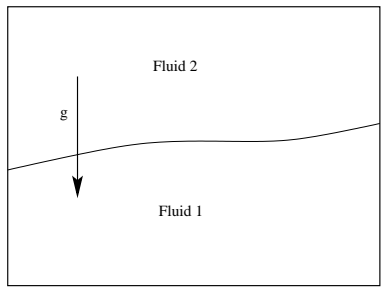,width=5cm}}
  \caption{Schematic representation of the two-fluid experiment with
    gravity.}
  \label{fig:grav}
\end{figure}
We consider two fluids subjected to a vertical gravity (see
Figure~\ref{fig:grav}). The lowest fluid is the heaviest. The domain
is $\Omega=(-2,2)\times(0,2)$. The equation of the steady state
interface is $x_2=1$. At $t=0$ the interface is perturbated, its
equation is $x_2=x_1/5 +1$. The experiment consists in observing
how the interface goes back to equilibrium (namely zero velocity,
hydrostatic pressure, interface $x_2=1$).  Here are the physical
parameters (in reduced units): $\rho_1=1$, $\rho_2=0.91$, $\eta_1=0.01$, $\eta_2=0.0091$,
$g=100$. In this test case, we neglect the surface tension effect
($\gamma=0$) and we assume pure slip on the wall (the GNBC will be
investigated in the next test case).

\subsubsection{Explicit treatment of the displacement of the interface}
\begin{figure}[htbp]
  \centering
  \centerline{\epsfig{file=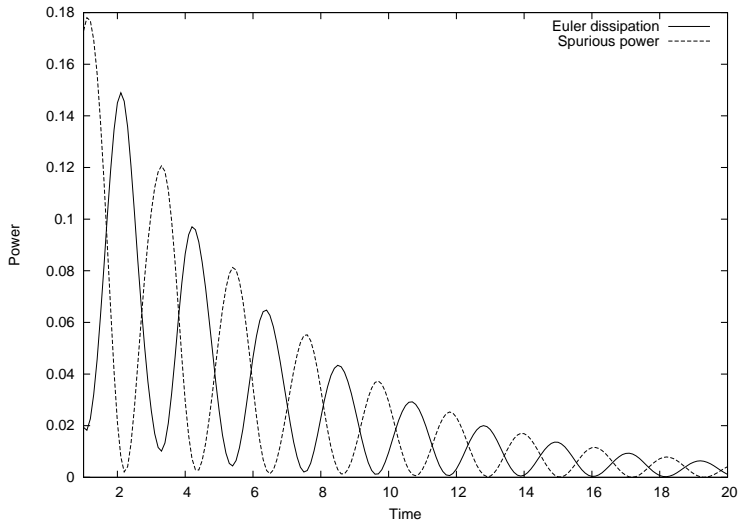,angle=-90,width=10cm}}
  \caption{Comparison of the spurious power term (\ref{eq:balance-expl})
    and the Euler dissipation (\ref{eq:euler}) with the explicit
    scheme (M1). The positiveness of the spurious power are in
    agreement with the results of Proposition~\ref{prop:balance-2} and can
    potentially induce unstabilities. The positiveness of the Euler
    dissipation have of course a stabilizing effect.}
  \label{fig:balance-expl}
\end{figure}

We first consider the explicit scheme (M1) with the body force
integrated on $\Om^{n+1}$ (namely $\Om^\star = \Om^{n+1}$).  In
Figure~\ref{fig:balance-expl}, we plot the quantity
\begin{equation}
  \label{eq:balance-expl}
\epsilon^n_{g,expl} = \Frac{K^{n+1} - K^n}{\delta t} + 
\Frac{W^{n+2} - W^{n+1}}{\delta t} + P_v^{n+1}
+ \Int{\Om^n}{}\Frac{\rho}{2\delta t}|\u^{n+1} - \u^n|^2,
\end{equation}
and the dissipation due to the Euler scheme
\begin{equation}
  \label{eq:euler}
  \epsilon^n_{euler}=\Int{\Om^n}{}\Frac{\rho}{2\delta t}|\u^{n+1} - \u^n|^2. 
\end{equation}
Two comments are in order. First, we observe that the quantity
$\epsilon^n_{g,expl}$ is indeed nonnegative which is an numerical
illustration of the result proved in
Proposition~\ref{prop:balance-2}. Second, we observe that the spurious
energy provided by the gravity term can be greater than the energy
dissipated by the Euler scheme. Other simulations (not reported here)
confirmed that both terms are of order $O(\delta t)$.

\begin{figure}[htbp]
  \centering
  \centerline{\epsfig{file=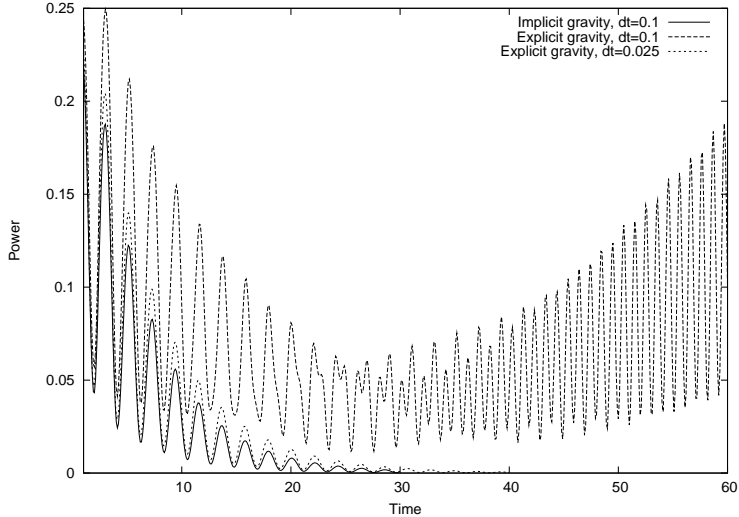,angle=-90,width=10cm}}
  \caption{Comparison of the viscous power term $P_v^n$ when the gravity
    is integrated on $\Om^\star=\Om^n$ or $\Om^{n+1}$ (referenced respectively as
    "explicit gravity" and "implicit gravity"). Although both choices
    have the same computational cost (we use scheme (M1)), we observe that the choice $\Om^\star=\Om^n$ may induce very severe unstabilities (here when $\delta t=0.1$).}
  \label{fig:explicit-grav-explicit}
\end{figure}

To illustrate the importance of the choice of the domain where the
body force is integrated, we have represented in Figure
\ref{fig:explicit-grav-explicit} the result obtained with $\Om^\star =
\Om^{n}$. It is worth noticing that the computional cost is the same
in both cases ($\Om^\star = \Om^{n}$ or $\Om^{n+1}$) since the
movement of the interface is treated explicitly in scheme (M1). With
$\delta t=0.025$ the result with $\Om^\star = \Om^{n}$ is similar to
the result with $\Om^\star = \Om^{n+1}$ and $\delta t=0.1$. But with
$\delta t=0.1$, we observe that the result may be dramatically
unstable when $\Om^\star = \Om^n$ (we indeed observe growing oscillations
of the free interface whereas the system is expected to go to its stationnary rest state).

\subsubsection{Implicit treatment of the displacement of the interface}
\begin{figure}[htbp]
  \centering
  \centerline{\epsfig{file=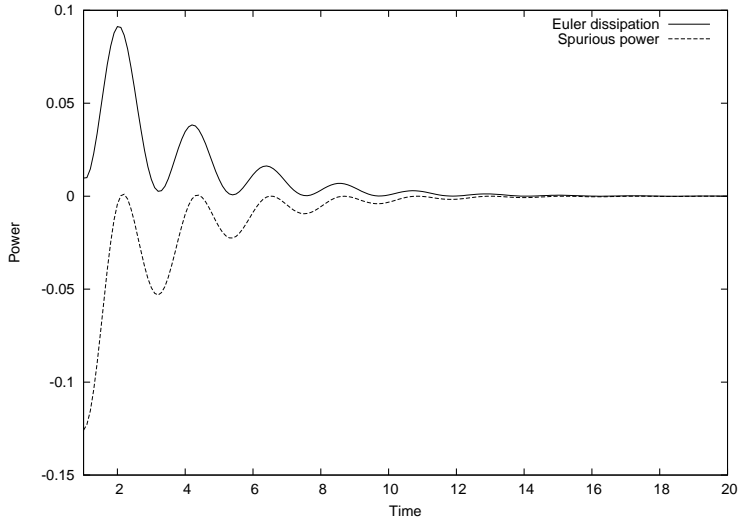,angle=-90,width=10cm}}
  \caption{Comparison of the spurious power term (\ref{eq:balance-impl})
    and the Euler dissipation (\ref{eq:euler}) with the (expensive)
    implicit scheme (M2). Contrarily to the result of Figure~\ref{fig:balance-expl}
    the balance is negative, which is in agreement with the results of
    Proposition~\ref{prop:balance-3} and ensure the stability of the
    scheme.}
  \label{fig:implicit_grav}
\end{figure}

In Figure~\ref{fig:implicit_grav}, we consider the fully implicit
scheme (M2) with $\Om^\star = \Om^{n+1}$ and we plot the quantity
\begin{equation}
  \label{eq:balance-impl}
  \epsilon^n_{g,impl} = \Frac{K^{n+1} - K^n}{\delta t} + 
  \Frac{W^{n+1} - W^{n}}{\delta t} + P_v^{n+1}
  + \Int{\Om^n}{}\Frac{\rho}{2\delta t}|\u^{n+1} - \u^n|^2,
\end{equation}
and the dissipation due to the Euler scheme. The results are in
agreement with Proposition~\ref{prop:balance-3}: the balance is now
negative which ensure the stability of the scheme. We nevertheless
recall that this interesting feature is obtained to the price of an
expensive scheme (since a nonlinear problem has to be solved at each
time step to achieve an implicit resolution of the interface
movement).

\subsubsection{Alternative schemes}

In this section we try to see how the theoretical results established
so far can lead to simple alternatives to the schemes (M1) and (M2).
The purpose is to devise a scheme whose cost is the same as (M1) and
which has dissipation properties similar to (M2).

Comparing equations (\ref{eq:gcl1_grav}) and (\ref{eq:gcl2_grav}) in
Lemma~\ref{lem:gcl_grav}, we observe that the last terms have opposite
signs: in one case it is dissipative, in the other case it brings a
spurious energy. It is therefore natural to try to combine these two
equations in order to decrease the amount of spurious energy appearing
in explicit schemes.  This simple observation leads to propose to
integrate a part of the gravity on $\Om^{n}$ and the other part on
$\Om^{n+1}$. In other words, we can take $\Om^\star=\Om^{n+1/2}$ in
(\ref{eq:weak-NS-disc}).

A second natural idea to circumvent the expensive implicit resolution
of the free interface movement is to use an explicit scheme with an
extrapolated interface velocity. More precisely, we propose to approximate~(\ref{eq:lapl-w-m2}) with:

Scheme (M3): displacement of the interface with an extrapolated velocity
\begin{equation}\label{eq:lapl-w-m2}
  \left\{
    \begin{array}{rcll}
      -\Delta w^n &= & 0,  & \mbox{ in } \Omega_i^n, i=1,2,\\
      w^n &= &\Frac{( 2\u^{n} - \u^{n-1} \circ \invalenm ) \cdot\n_{\Sigma^{n}}  }{n_{\Sigma^{n}}^{3}},  
      &\mbox{ on } \Sigma^{n}, \\
      \dps{\frac{\partial w^n}{\partial \n}} &=& 0, & \mbox{ on } \partial \Omega,
    \end{array}
  \right.
\end{equation}

\begin{figure}[htbp]
  \centerline{\epsfig{file=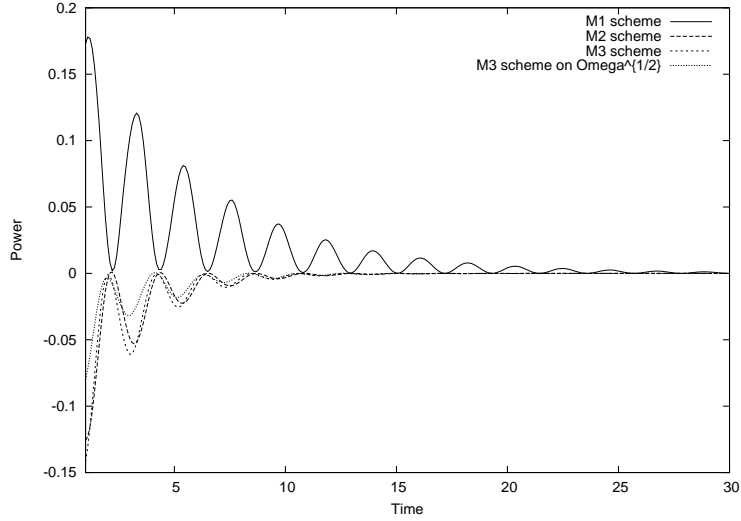,angle=-90,width=10cm}}
  \caption{Spurious power term (\ref{eq:balance-expl}) for scheme (M1)
    (same result as in Figure~\ref{fig:balance-expl}) and the term
    (\ref{eq:balance-impl}) for the implicit scheme (M2), or the
    explicit scheme (M3) with an extrapolated interface (with for all
    these computations $\Om^\star=\Om^{n+1}$). We observe
    that when the gravity is integrated on $\Om^\star=\Om^{n+1/2}$ (with
    the scheme (M3)), the scheme
    is still stable (negative energy balance) but also less dissipative.}
  \label{fig:extrap_grav_1_2}
\end{figure}

We have tested these two simple ideas in the above experiment (two
fluids submitted to gravity).  The results are reported in
Figure~\ref{fig:extrap_grav_1_2}. First, we set $\Om^\star=\Om^{n+1}$,
and we compare the scheme (M2) with the scheme (M3). We observe that they have almost the same dissipation
property. Thus, (M3) has (almost) the same stability properties as the
expensive (M2) scheme and the same computational cost as the "cheap"
(M1) scheme.  Second, using again scheme (M3), we notice that integrating the gravity on
$\Om^\star=\Om^{n+1/2}$ decreases the artificial dissipation of the
scheme. 

Therefore, it seems that a good compromise between efficiency, stability
and artificial dissipation is to use scheme (M3) with $\Om^\star=\Om^{n+1/2}$.

Of course, these properties still have to be assessed in more complex
situations. The only purpose of this section was to show the potential
interest of the energy studies presented above to devise new schemes.

\subsection{Moving contact line problems and GNBC}\label{sec:num-gnbc}

\begin{figure}[htbp]
  \centering
  \centerline{\epsfig{file=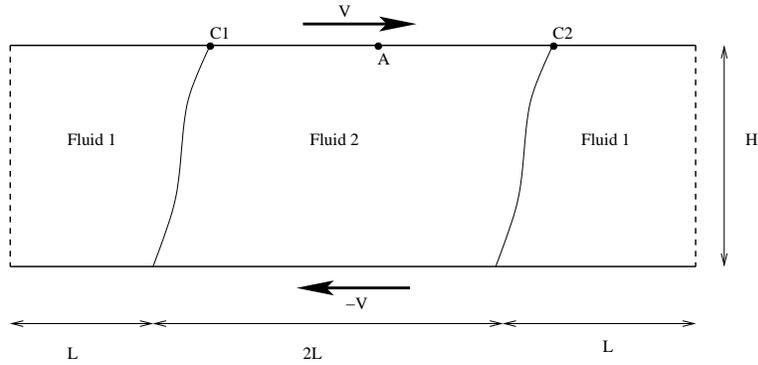,width=10cm}}
  \caption{Schematic representation of the two-fluid periodic Couette simulation.}
  \label{fig:couette-bif}
\end{figure}

To illustrate how the GNBC can handle the moving contact line, we
present the results of two benchmarks proposed in
\cite{qian-wang-sheng-06} involving a Couette flow for two fluids.
The geometry is 2D and represented in
Figure~\ref{fig:couette-bif}. 
The domain is periodic along $x$: the
dashed lines, located on $x=0$ and $x=4L$ represent the periodic
boundaries.  The walls are defined by $y=0$ (bottom) and $y=H$
(top). A velocity $V \uu{e}_x$ (resp.  $-V\uu{e}_x$) is imposed on the
top (resp. on the bottom) of the domain.  For the first test case (the
``symmetric'' one), we took the following values of the parameters
from \cite{qian-wang-sheng-06} (given in reduced units): $H=13.6$,
$L=27.2$, $V=0.25$, $\rho_1=\rho_2=0.81$, $\eta_1=\eta_2=1.95$,
$\gamma=5.5$, $\beta_1=\beta_2=1.5$ and $\theta_s=\pi/2$.  The
parameters for the second test case (the ``asymmetric'' one) are the
same except $V=0.20$, $\beta_2=0.591$ and $\theta_s$ such that
$\cos\theta_s\approx 0.38$. In both cases, of course, there is no
gravity: $g=0$.

\begin{figure}[htbp]
  \centering
  \centerline{\epsfig{file=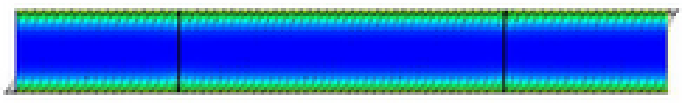,width=14cm}}
  \centerline{\epsfig{file=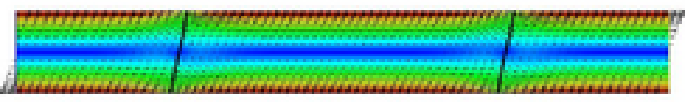,width=14cm}}
  \centerline{\epsfig{file=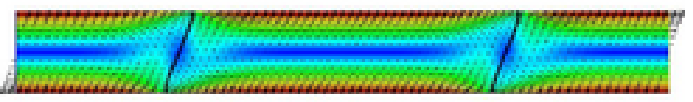,width=14cm}}
  \centerline{\epsfig{file=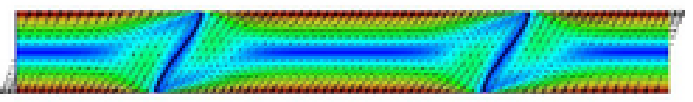,width=14cm}}
  \caption{Velocity field at t=1, 10, 20, 100 in the symmetric case. The color
    represents the velocity magnitude.}
  \label{fig:vel}
\end{figure}

\begin{figure}[htbp]
  \centering
  \centerline{\epsfig{file=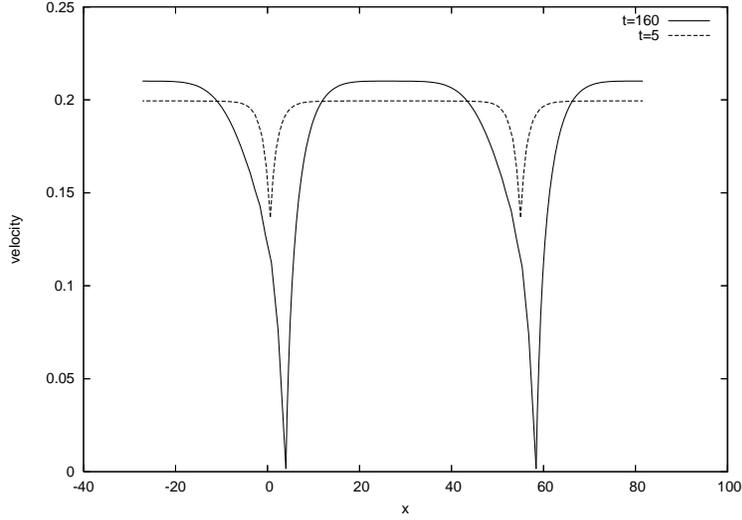,angle=-90,width=10cm}}
  \caption{Velocity on the top wall \emph{versus} $x$ for the symmetric
    case.  At $t=5$, we are still in the transient phase, the velocity of
    the points on the contact line is non-zero.  At $t=160$, the points on
    the contact line are almost fixed (which means they indeed move with
    respect to the wall), whereas the points on the boundary far from
    the contact line move with a velocity of magnitude about $0.21$
    ($0.25$ would have meant a complete adherence on the wall).}
  \label{fig:veltop}
\end{figure}

\begin{figure}[htbp]
  \centering
  \centerline{\epsfig{file=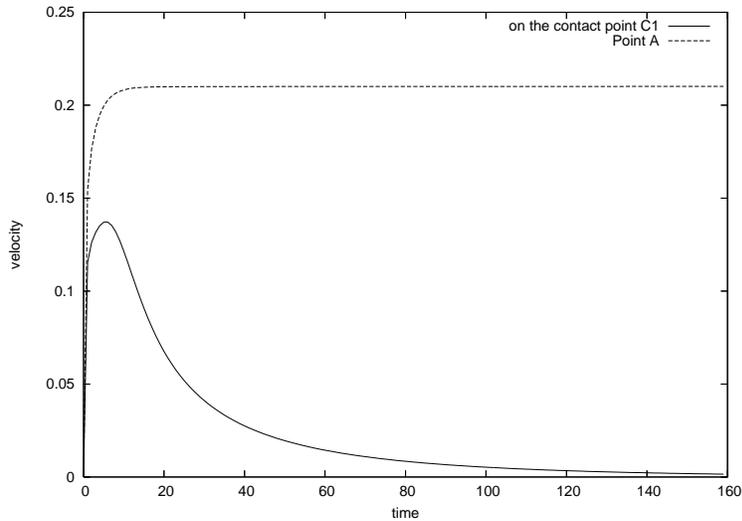,angle=-90,width=10cm}}
  \caption{Velocity of the point C1 on the contact line 
    and of the point A on the boundary (see Figure~\ref{fig:couette-bif})
    \emph{versus} time for the symmetric case. The velocity of the point
    on the contact line tends to zero (pure slip) whereas the velocity
    of the point A tends to about $0.21$.}
  \label{fig:veltime}
\end{figure}

\begin{figure}[htbp]
  \centering
  \centerline{\epsfig{file=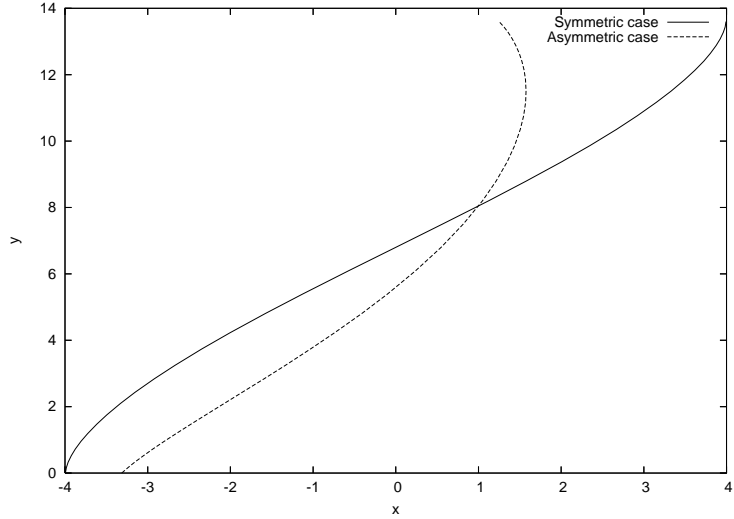,angle=-90,width=10cm}}
  \caption{Stationnary interface profiles in the symmetric and asymmetric cases.}
  \label{fig:interf}
\end{figure}

At $t=0$, the interfaces separating the two fluids are straight and
vertical. After a while the interfaces reach a steady state position.
Let us emphasize that this behaviour is a direct result of the GNBC
conditions: the interfaces would of course not converge to steady
curves if we had imposed $u_x=V$ on the top and $u_x=-V$ on the bottom
(no-slip). On the other hand, such a result cannot be obtained with
pure slip boundary conditions.  Figure~\ref{fig:vel} shows the
velocity field at $t=1, 10, 20, 100$ in the symmetric case.  The color
represents the magnitude of the velocity.  At $t=100$ the stationnary
state is almost reached: the blue zone surrounding the interfaces
shows that they are fixed (which means they indeed sleep with respect
to the wall), whereas the remaining part of the wall is red, which
corresponds to a fluid adherence on the wall.
More quantitatively, Figure~\ref{fig:veltop} shows the velocity on the
top wall at $t=5$ (the contact points are still moving) and at $t=160$
(the contact points are fixed).  Figure~\ref{fig:veltime} shows the
evolution in time of the velocity of a point on the contact line
(which tends to zero) and of the velocity of a point on the wall far
from the contact line (which tends to about $0.21$, whereas a total
adherence would correspond to $0.25$).  The stationary interfaces in
the symmetric and asymmetric cases are represented on
Figure~\ref{fig:interf}.  These results are in very good agreement
with those presented in~\cite{qian-wang-sheng-06}, which were obtained
either by a continuum phase-field formulation, or by molecular
dynamics simulations.


\begin{figure}[htbp]
  \centering
  \centerline{\epsfig{file=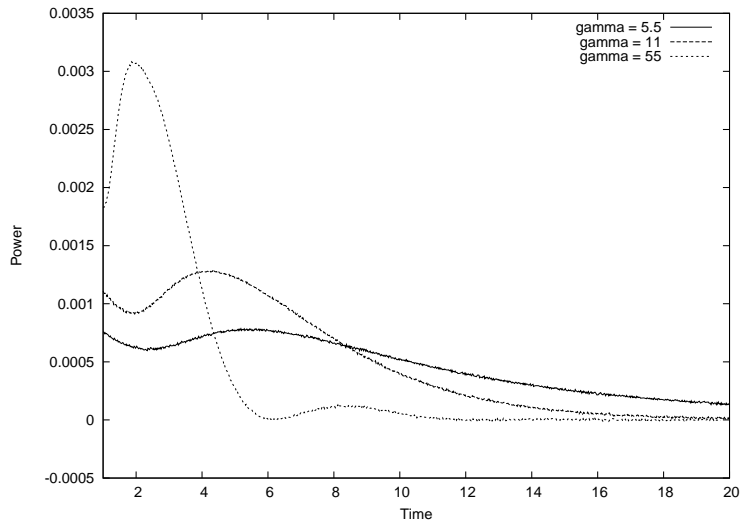,angle=-90,width=10cm}}
  \caption{Effect of surface tension on the energy balance
    (\ref{eq:balance_surf_gnbc}) when the interface displacement is treated
    explicitly. Observe that this quantity is indeed positive and
    increases with the surface tension coefficient $\gamma$ (at least
    in the dynamic part of the simulation when the interface is
    still significantly moving). This confirms the results of
    Proposition~\ref{prop:balance-2}.}
  \label{fig:energie-surf}
\end{figure}

To illustrate the result established in
Proposition~\ref{prop:balance-2}, we performed the above experiment
with three values of surface tension $\gamma=5.5, 11, 55$ and we plot
in Figure~\ref{fig:energie-surf}, the quantity:
\begin{equation}
  \begin{array}[t]{l}
    \Frac{K^{n+1} - K^n}{\delta t} + 
    \Frac{W^{n+2} - W^{n+1}}{\delta t} + P_v^{n+1}
    + \gamma\Frac{|\Sigma^{n+2}| - |\Sigma^{n+1}|}{\dt} 
    + \beta(\u^{n+1}-\u_b,\u^{n+1})_{\partial\Om}
    \\[9pt]   ~~~~~~~~+ \Int{\Om^n}{}\Frac{\rho}{2\delta t}|\u^{n+1} - \u^n|^2 
    - \gamma \int_{\partial\Sigma^{n+1}} \cos(\theta_s) \t_{\partial\Om}\cdot \u^{n+1}\,dl_{\partial\Sigma}.
  \end{array}
  \label{eq:balance_surf_gnbc}
\end{equation} 
According to Proposition~\ref{prop:balance-2}, this quantity is equal
to $\epsilon_{\gamma,exp}^{n+1}$ (defined in
\eqref{eq:epsilon-gamma}). The results presented in
Figure~\ref{fig:energie-surf} confirm that the surface tension indeed
generates spurious power (the balance is positive), which
increases with the surface tension coefficient $\gamma$.

\section{Conclusion}

We have proposed a formulation for the Generalized Navier Boundary
Condition which allows to compute the moving contact line at the
interface between two fluids. This formulation is in particular
well-suited to an energy stability analysis. We have shown that an
explicit treatment of the free interface displacement introduces a spurious power
in presence of gravity and surface tension. The study of the surface
tension terms required an extension of the Geometric Conservation Law
concept to surface integrals.  Several stability results have been
established and illustrated with numerical experiments. Understanding the
spurious power induced by numerical schemes may help to devise more
stable algorithms. In particular, we have proposed a simple scheme
which offers a good compromise between efficiency, stability and
artificial diffusion.

Extension to second order schemes and generalization to more general
cases of the ``Surface Geometric Conservation Law''
(Lemma~\ref{lem:gcl_surf}) could be investigated in future works.


\end{document}